\documentclass[11pt]{article}

\usepackage{color}
\usepackage{latexsym}
\usepackage{amssymb}
\usepackage{amsmath,amsfonts,theorem,euscript,array,enumerate,mathrsfs}
\usepackage{graphicx}

\newtheorem{Theorem}{Theorem}[section]
\newtheorem{Definition}[Theorem]{Definition}
\newtheorem{Proposition}[Theorem]{Proposition}

\newtheorem{Remark}[Theorem]{Remark}

\numberwithin{equation}{section}

\def\esssup_#1{\underset{#1}{\mathrm{ess\,sup\, }}}
\def\essinf_#1{\underset{#1}{\mathrm{ess\,inf\, }}}

\def \N{\mathbb{N}}
\def \R{\mathbb{R}}

\def \E{\mathbb{E}}
\def \F{\mathbb{F}}

\def \H{\mathbb{H}}

\def \P{\mathbb{P}}

\def \D{\mathbb{D}}
\def \S{\mathbb{S}}

\def \X{\mathbb{X}}

\def \Bc{{\cal B}}

\def \Dc{{\cal D}}

\def \Fc{{\cal F}}

\def \Mc{{\cal M}}

\def \Sc{{\cal S}}

\def \Uc{{\cal U}}

\def \Xc{{\cal X}}

\def \eps{\varepsilon}

\def \ep{\hbox{ }\hfill$\Box$}

\allowdisplaybreaks

\begin{document}

\title{\small{\textbf{FUNCTIONAL IT\^O VERSUS BANACH SPACE STOCHASTIC CALCULUS\\ AND STRICT SOLUTIONS OF SEMILINEAR PATH-DEPENDENT EQUATIONS}}}

\author{
\textsc{Andrea Cosso}\thanks{Laboratoire de Probabilit\'es et Mod\`eles Al\'eatoires, CNRS, UMR 7599, Universit\'e Paris Diderot, France. e-mail: \sf andrea.cosso at math.univ-paris-diderot.fr}
\qquad\quad
\textsc{Francesco Russo}\thanks{Unit\'e de Math\'ematiques Appliqu\'ees, 
ENSTA ParisTech, Universit\'e Paris-Saclay,
828, boulevard des Mar\'echaux, F-91120 Palaiseau, France. e-mail: \sf francesco.russo at ensta-paristech.fr}
}

\date{April 24, 2015}

\maketitle

\begin{abstract}
\noindent Functional It\^o calculus was introduced
in order to expand a functional $F(t, X_{\cdot+t}, X_t)$ depending on time $t$, past and present values of the process $X$.
Another possibility to expand $F(t, X_{\cdot+t}, X_t)$ consists in considering 
the path $X_{\cdot+t}=\{X_{x+t},\,x\in[-T,0]\}$ as an element of the Banach space of continuous functions on $C([-T,0])$
and to use Banach space stochastic calculus.
The aim of  this paper is threefold. 
1) To reformulate functional It\^o calculus, separating time and past, making use of the regularization procedures which matches more naturally the notion of horizontal derivative which is one of
the tools of that calculus. 
2) To exploit this reformulation in order to discuss the (not obvious) relation between the functional and the Banach space approaches.
3) To study existence and uniqueness of smooth solutions to path-dependent partial differential equations which naturally arise in the study of functional It\^o calculus. More precisely, we study a path-dependent equation of Kolmogorov type which is related to the window process of the solution to an It\^o stochastic differential equation with path-dependent coefficients. We also study a semilinear version of that equation. 

\vspace{4mm}

\noindent {\bf Keywords:} functional It\^o calculus; 
 Banach space valued stochastic calculus;  
path-dependent partial differential equation; strict solutions; 
calculus via regularization. 

\vspace{4mm}

\noindent {\bf AMS 2010 subject classifications:}  
 60H05; 60H10; 60H30; 35A09; 35K10.
\end{abstract}

\section{Introduction}

Recently, a new branch of  stochastic calculus has appeared, known as \emph{functional It\^o calculus}, which results to be an extension of classical It\^o calculus to 
functionals depending on the all path of a stochastic process and not only on its current value, see Dupire \cite{dupire}, Cont and Fourni\'e 
\cite{contfournie10,contfournie,contfournie13}.
On the other hand,  C. Di Girolami, the second named author
 and more recently
G. Fabbri, have introduced a path-dependent type calculus having similar objectives in a series of papers (\cite{DGR, DGRnote, 
digirrusso12, DGR2, digirfabbrirusso13}), which is a stochastic calculus
 for processes taking values in
a separable Banach space $B$ (which includes the case $B = C([-T,0])$).

\vspace{1mm}

\noindent The aim of  the paper is threefold. 
\begin{enumerate}
\item First, to reformulate functional It\^o calculus, separating time and past, making use of the regularization procedures which matches more naturally the notion of horizontal derivative which is one of
the tools of that calculus. 
\item Exploiting this reformulation of functional It\^o calculus in order to discuss the relation between the functional and the Banach space approaches. In particular, to investigate the (not obvious) nature of the so-called horizontal derivative. 
\item Study path-dependent Partial Differential Equations (PDEs) of Kolmogorov type which are associated to window processes of solutions to It\^o stochastic differential equations with path-dependent drift and diffusion coefficients.
\end{enumerate}
The paper is separated into two sections/parts: items 1. and 2. above are in the first part, while item 3. is studied in the second part. Let us now describe more in detail these two sections.

\vspace{1mm}

In the first part of the paper, we revisit functional It\^o calculus by means of stochastic calculus via regularization. We recall that Cont and Fourni\'e  \cite{contfournie10,contfournie,contfournie13} 
developed functional It\^o calculus and derived a functional It\^o formula 
using discretization techniques of F\"ollmer \cite{follmer} type, instead of regularization techniques. One of the main issues of functional It\^o's calculus is the definition of the functional (or pathwise) derivatives, i.e., the
 horizontal derivative (calling in only the past values of the trajectory) and the vertical derivative (calling in only the present value of the trajectory). 
As mentioned above, the idea of regularization makes the approach more natural since it matches better the notion
of derivative (for instance horizontal). On the other hand we have decided to keep separated
``time'' and ``past'', the sum of the time derivative plus ``our'' horizontal derivative would
give Dupire's horizontal derivative.
In \cite{contfournie10}, it is essential to consider functionals defined on the space of c\`adl\`ag trajectories, namely $\D([-T,0])$, since the definition of functional derivatives necessitates of discontinuous paths. Therefore, if a functional is defined only on the space of continuous trajectories (because, e.g., it depends on the paths of a continuous process as the Brownian motion), we have to extend it anyway to the space of c\`adl\`ag trajectories, even though, in general, there is no 
a unique way to extend it. 
In contrast with this approach, our point of view is to consider a state space which get stuck as much as possible to the ``natural space'' of continuous functions $C([-T,0])$. In the classical literature of functional
dependent stochastic differential equations (as delay equations, see for instance \cite{cho}), 
the process takes values into $L^2 ([-T,0]) \times \R$, the first component staying
for the past of the trajectory and the second one for the present.
This space is isomorphic to $L^2([-T,0], d\mu)$ where $\mu$ is the sum of Lebesgue measure
and Dirac measure at zero. This space contains strictly $C([-T,0])$ as a dense subspace.
Our idea was to consider as state space $\mathscr C([-T,0])$, that is a space which contains $C([-T,0])$ as a dense subset and
is isomorphic to a product space which allows to separate past and present. On $\mathscr C([-T,0])$ we define the functional derivatives. 
$\mathscr C([-T,0])$ is the space of bounded trajectories on $[-T,0]$, continuous on $[-T,0[$ and with possibly a jump at $0$. We endow $\mathscr C([-T,0])$ with a topology such that $C([-T,0])$ is dense in $\mathscr C([-T,0])$ with respect to this topology. Therefore, any functional $\Uc\colon[0,T]\times C([-T,0])\rightarrow\R$, continuous with respect to the topology of $\mathscr C([-T,0])$, admits at most a unique extension to $\mathscr C([-T,0])$, denoted $u\colon[0,T]\times\mathscr C([-T,0])\rightarrow\R$. We present some significant functionals for which a continuous extension exists. Then, we develop the functional It\^o calculus for $u\colon[0,T]\times\mathscr C([-T,0])\rightarrow\R$.

Notice that we use a slightly different notation with respect to \cite{contfournie10}. In particular, instead of a map $\Uc\colon [0,T]\times C([-T,0])\rightarrow\R$, in \cite{contfournie10} a family of maps $F=(F_t)_{t\in[0,T]}$, with $F_t\colon C([0,t])\rightarrow\R$, is considered. However, we can always move from one formulation to the other. Indeed, given $F=(F_t)_{t\in[0,T]}$, where each $F_t\colon C([0,t])\rightarrow\R$, we can define $\Uc\colon [0,T]\times C([-T,0])\rightarrow\R$ as follows:
\[
\Uc(t,\eta) \ := \ F_t(\eta(\cdot+T)|_{[0,t]}), \qquad (t,\eta)\in[0,T]\times C([-T,0]).
\]
Vice-versa, let $\Uc\colon [0,T]\times C([-T,0])\rightarrow\R$ and define $F=(F_t)_{t\in[0,T]}$ as
\begin{equation}
\label{E:F=u}
F_t(\tilde\eta) \ := \ \Uc(t,\eta), \qquad (t,\tilde\eta)\in[0,T]\times C([0,t]),
\end{equation}
where $\eta$ is the element of $C([-T,0])$ obtained from $\tilde\eta$ first translating $\tilde\eta$ on the interval $[-t,0]$, then extending it in a constant way up to $-T$, namely $\eta(x) := \tilde\eta(x+t)1_{[-t,0]}(x) + \tilde\eta(-t)1_{[-T,-t)}(x)$, for any $x\in[-T,0]$. Observe that, in principle, the map $\Uc$ contains more information than $F$, since in \eqref{E:F=u} we do not take into account the values of $\Uc$ at $(t,\eta)\in[0,T]\times C([-T,0])$ with $\eta$ not constant on the interval $[-T,-t]$. Despite this, the equivalence between the two notations is guaranteed by the fact that, as it will be clear later, when we consider the composition of $\Uc$ with a stochastic process, this extra information plays no role. Our formulation has two advantages. Firstly, we can work with a single map instead of a family of maps. In addition, the time variable and the path have two distinct roles in our setting, as for the time variable and the space variable in the classical It\^o calculus. This, in particular, allows us to define the horizontal derivative independently of the time derivative, so that, the horizontal derivative defined in \cite{contfournie10} corresponds to the sum of our horizontal derivative and of the time derivative. We mention that an alternative approach to functional derivatives was introduced in \cite{buckdahn_ma_zhang13}.

After this reformulation of functional It\^o calculus, we can now investigate the relation between functional It\^o calculus and Banach space valued stochastic calculus (via regularization), for the case of window processes. This latter and brand new branch of stochastic calculus and stochastic analysis has been recently conceived, deeply studied, and developed in many directions in  \cite{DGR2, digirrusso12, DGRnote},
\cite{digirfabbrirusso13} and for more details \cite{DGR}. 
For the particular case of window processes, we also refer to Theorem 6.3 and Section 7.2 in \cite{digirfabbrirusso13}. In the present paper, we prove formulae which allow to express functional derivatives in terms of differential operators arising in the Banach space valued stochastic calculus via regularization. In particular, while the identification of the vertical derivative is rather expected, we found a not obvious relation between the horizontal derivative and second order derivative operators of Banach space valued stochastic calculus.

Dupire \cite{dupire} introduced also the concept of \emph{path-dependent partial differential equation}, to which the second part of the present paper is devoted. Di Girolami and the second named author, in Chapter 9
of \cite{DGR}, considered
a similar equation in the framework of Banach space valued calculus, 
for which we refer also to \cite{flandoli_zanco13}. 
We also drive the attention to the recent contribution of
  \cite{leao_ohashi_simas14}.

In the last part of the paper we focus on path-dependent semilinear Kolmogorov equations driven by a path dependent generator associated with a delay equation, for which we provide a definition of strict solution (namely smooth solution; we prefer to use the term ``strict'' instead of ``classical'' because all the theory of path-dependent partial differential equations is very recent). We prove a uniqueness result for this kind of solution, by means of probabilistic methods based on the theory of Backward Stochastic Differential Equations (BSDEs). More precisely, we show that, if a strict solution exists, then it can be expressed through the  solution of
 a certain backward stochastic differential equation. Therefore, from the uniqueness of the BSDE it follows that there exists at most one strict solution. Then,  we also prove an existence result for strict solutions. 
In the companion paper \cite{cosso_russo15b}, for the same type of equations, we introduce a more
general notion of solution, that we have denominated {\it strong-viscosity solution}, for which we provide existence and uniqueness results.

The  paper is organized as follows. In section \ref{S:Functional} we develop functional It\^o calculus via regularization: after a brief introduction on finite dimensional stochastic calculus via regularization in subsection \ref{SubS:Background}, we introduce and study the space $\mathscr C([-T,0])$ in subsection \ref{SubS:Preliminaries}; then, we define the functional derivatives and we prove the functional It\^o formula in subsection \ref{SubS:PathwiseDerivatives}; in subsection \ref{SubS:Comparison}, instead, we discuss the relation between functional It\^o calculus via regularization and Banach space valued stochastic calculus via regularization for window processes. Finally, in Section \ref{S3} we study strict solutions to path-dependent PDEs.

\section{Functional It\^o calculus: a regularization approach}
\label{S:Functional}

\subsection{Background: finite dimensional calculus via regularization}
\label{SubS:Background}

The theory of stochastic calculus via regularization has been developed in several papers, starting from \cite{russovallois91, russovallois93}. We recall below only the results used in the present paper, and we refer to \cite{russovallois07} for a survey on the subject. We emphasize that integrands are allowed to be anticipating. Moreover, the integration theory and calculus appears to be close to a pure pathwise approach even though there is still a probability space behind.

\vspace{3mm}

Fix a probability space $(\Omega,\Fc,\P)$ and $T\in]0,\infty[$. Let $\F=(\Fc_t)_{t\in[0,T]}$ denote a filtration satisfying the usual conditions. Let $X=(X_t)_{t\in[0,T]}$ (resp. $Y=(Y_t)_{t\in[0,T]}$) be a real continuous (resp. $\P$-a.s. integrable) process. Every real continuous process $X = (X_t)_{t\in[0,T]}$ is naturally extended to all $t\in\R$ setting $X_t = X_0$, $t\leq0$, and $X_t=X_T$, $t\geq T$. We also define a $C([-T,0])$-valued process $\X=(\X_t)_{t\in\R}$, called the {\bf window process} associated with $X$, defined by
\begin{equation}
\label{WindowProcess}
\X_t := \{X_{t+x},\,x\in[-T,0]\}, \qquad t\in\R.
\end{equation}

\begin{Definition}
\label{D:DPFI}
Suppose that, for every $t\in[0,T]$,  the following limit
\begin{equation}
\label{DPFI}
\int_{0}^{t}Y_sd^-X_s \ := \ \lim_{\eps\rightarrow0^+}\int_{0}^{t} Y_s\frac{X_{s+\eps}-X_s}{\eps}ds,
\end{equation}
exists in probability. If the obtained  random function admits a continuous modification, that process is denoted by  $\int_0^\cdot Yd^-X$ and called \textbf{forward integral of $Y$ with respect to $X$}.
\end{Definition}

\begin{Definition}
A family of processes $(H_t^{(\eps)})_{t\in[0,T]}$ is said to converge to $(H_t)_{t\in[0,T]}$ in the \textbf{ucp sense}, if $\sup_{0\leq t\leq T}|H_t^{(\eps)}-H_t|$ goes to $0$ in probability, as $\eps\rightarrow0^+$.
\end{Definition}

\begin{Proposition}
Suppose that the limit \eqref{DPFI} exists in the ucp sense. Then, the forward integral $\int_0^\cdot Yd^-X$ of $Y$ with respect to $X$ exists.
\end{Proposition}

Let us introduce the concept of covariation, which is a crucial notion in stochastic calculus via regularization. Let us suppose that
$X, Y$ are continuous processes.

\begin{Definition}
\label{D33}
The \textbf{covariation of $X$ and $Y$} is defined by
\[						
\left[X,Y\right]_{t} \ = \ \left[Y,X\right]_{t}  \ = \ \lim_{\eps\rightarrow 0^{+}}
\frac{1}{\eps} \int_{0}^{t} (X_{s+\eps}-X_{s})(Y_{s+\eps}-Y_{s})ds, \qquad t \in [0,T],
\]
if the limit exists in probability for every  $t \in [0,T]$, provided that the limiting random function  admits a continuous version $($this is the case if the limit holds in the ucp sense$)$. If $X=Y,$ $X$ is said to be a \textbf{finite quadratic variation process} and we set $[X]:=[X,X]$.
\end{Definition}

The forward integral and the covariation generalize the classical It\^o integral and covariation for semimartingales. In particular, we have the following result.

\begin{Proposition}
\label{Pproperties}
The following properties hold.
\begin{enumerate}
\item[\textup{(i)}] Let $S^1,S^2$ be continuous $\F$-semimartingales. Then, $[S^1,S^2]$ is the classical bracket $[S^1,S^2]=\langle M^1,M^2\rangle$, where $M^1$ $($resp. $M^2$$)$ is the local martingale part of $S^1$ $($resp. $S^2$$)$.
\item[\textup{(ii)}] Let $V$ be a continuous bounded variation process and $Y$ be a c\`{a}dl\`{a}g process $($or vice-versa$)$$;$ then $[V] =[Y,V]= 0$. Moreover $\int_0^\cdot Y d^-V=\int_0^\cdot Y dV $,
is the {\bf Lebesgue-Stieltjes integral}.
\item[\textup{(iii)}] 
If $M$ is an  $\F$-local martingale and $Y$ is an $\F$-progressively measurable c\`{a}dl\`{a}g process, then $\int_0^\cdot Yd^- M$ exists and equals the It\^o integral  $\int_0^\cdot YdM$.
\end{enumerate}
\end{Proposition}
\textbf{Proof.}
Assertion (i) follows from Corollary 2 and Proposition 9 in \cite{russovallois07}. Concerning (ii) we refer to item 7) of Proposition 1 in \cite{russovallois07}. (iii) follows from Proposition 6 in \cite{russovallois07}.
\ep

\vspace{3mm}

We end this crash introduction to finite dimensional stochastic calculus via regularization presenting one of its cornerstones: It\^o's formula. It is a  well-known result in the theory of semimartingales, but it also  extends to the framework of  finite quadratic variation processes. For a proof we refer to Theorem 2.1 of \cite{russovallois95}.

\begin{Theorem}
\label{ITOFQV}
Let $F:[0,T]\times \R\longrightarrow \R$ be of class $ C^{1,2}\left( [0,T]\times \R \right)$ and $X=(X_t)_{t\in[0,T]}$ be
a real continuous finite quadratic variation process. Then, the following \textbf{It\^o formula} holds:
\begin{align}
F(t,X_{t}) \ &= \ F(0,X_{0}) + \int_{0}^{t}\partial_t  F(s,X_s)ds + \int_{0}^{t} \partial_{x} F(s,X_s)d^-X_s + \frac{1}{2}\int_{0}^{t} \partial^2_{x\, x} F(s,X_s)d[X]_s, \label{E:ITOFQV}
\end{align}
for all $0\leq t\leq T$.
\end{Theorem}

\subsubsection{The deterministic calculus via regularization}

\label{S211}

A useful particular case of finite dimensional stochastic calculus via regularization arises when $\Omega$ is a singleton, i.e., when the calculus becomes deterministic. In addition, in this deterministic framework we will make use of the \emph{definite integral} on an interval $[a,b]$, where $a<b$ are two real numbers. Typically, we will consider $a=-T$ or $a=-t$ and $b=0$.

\vspace{3mm}

We start with two conventions. By default, every bounded variation function $f\colon[a,b]\rightarrow\R$ will be considered as c\`{a}dl\`{a}g. Moreover, given a function $f\colon[a,b]\rightarrow\R$, we will consider the following two extensions of $f$ to the entire real line:
\[
f_J(x) :=
\begin{cases}
0, & x>b, \\
f(x), & x\in[a,b], \\
f(a), & x<a,
\end{cases} \qquad\qquad
f_{\overline J}(x) :=
\begin{cases}
f(b), & x>b, \\
f(x), & x\in[a,b], \\
0, & x<a,
\end{cases}
\]
where $J:=\,]a,b]$ and $\overline J=[a,b]$.

\begin{Definition}
\label{D:DeterministicIntegral}
Let $f,g\colon[a,b]\rightarrow\R$ be c\`{a}dl\`{a}g functions.\\
\textup{(i)} Suppose that the following limit
\[
\int_{[a,b]}g(s)d^-f(s) \ := \ \lim_{\eps\rightarrow0^+}\int_\R g_J(s)\frac{f_{\overline J}(s+\eps)-f_{\overline J}(s)}{\eps}ds,
\]
exists and it is finite. Then, the obtained quantity is denoted by $\int_{[a,b]} gd^-f$ and called \textbf{$($determin\-istic, definite$)$ forward integral of $g$ with respect to $f$ $($on $[a,b]$$)$}.\\
\textup{(ii)} Suppose that the following limit
\[
\int_{[a,b]}g(s)d^+f(s) \ := \ \lim_{\eps\rightarrow0^+}\int_\R g_J(s)\frac{f_{\overline J}(s)-f_{\overline J}(s-\eps)}{\eps}ds,
\]
exists and it is finite. Then, the obtained quantity is denoted by $\int_{[a,b]} gd^+f$ and called \textbf{$($determin\-istic, definite$)$ backward integral of $g$ with respect to $f$ $($on $[a,b]$$)$}.
\end{Definition}

The notation concerning this integral is justified by the fact that when the integrator
$f$ has bounded variation then previous integrals are Lebesgue-Stieltjes
integrals on $[a,b]$, as stated in the following proposition, whose simple proof is not reported.
\begin{Proposition}
Suppose $f\colon[a,b]\rightarrow\R$ with bounded variation and $g\colon[a,b]\rightarrow\R$ c\`{a}dl\`{a}g. Then, we have
\begin{align*}
\int_{[a,b]} g(s) d^-f(s) \ &= \ \int_{[a,b]} g(s^-) df(s) \ := \ g(a)f(a) + \int_{]a,b]} g(s^-) df(s), \\
\int_{[a,b]} g(s) d^+f(s) \ &= \ \int_{[a,b]} g(s) df(s) \ := \ g(a)f(a) + \int_{]a,b]} g(s) df(s).
\end{align*}
\end{Proposition}
Let us now introduce the deterministic covariation.

\begin{Definition}
\label{D:DeterministicQuadrVar}
Let $f,g\colon[a,b]\rightarrow\R$ be continuous functions and suppose that $0\in[a,b]$. The \textbf{$($deterministic$)$ covariation of $f$ and $g$ $($on $[a,b]$$)$} is defined by
\[						
\left[f,g\right](x) \ = \ \left[g,f\right](x)  \ = \ \lim_{\eps\rightarrow 0^{+}}
\frac{1}{\eps} \int_{0}^x (f(s+\eps)-f(s))(g(s+\eps)-g(s))ds, \qquad x\in[a,b],
\]
if the limit exists and it is finite for every  $x\in[a,b]$. If $f=g$, we set $[f]:=[f,f]$ and it is called \textbf{(deterministic) quadratic variation of
 $f$ $($on $[a,b]$$)$}.
\end{Definition}

We notice that in Definition \ref{D:DeterministicQuadrVar} the quadratic variation $[f]$ is continuous on $[a,b]$, since $f$ is a continuous function.

\begin{Remark}
{\rm
Notice that if $f$ is a fixed Brownian path and $g(s)=\varphi(s,f(s))$, with $\varphi\in C^1([0,1]\times \R)$ (for simplicity, we take $[a,b]=[0,1]$), then $\int_{[0,1]} g(s) d^- f(s)$ exists for almost all (with respect to the Wiener measure on $C([0,1])$) Brownian path $f$. This latter result can be shown using Theorem 2.1 in \cite{gradinaru_nourdin03} (which implies that the deterministic bracket, introduced in Definition \ref{D:DeterministicQuadrVar} below, exists, for almost all Brownian paths $f$, and $[f](s)=s$) and then applying It\^o's formula in Theorem \ref{ITOFQV} above, with $\P$ given by the Dirac delta at a Brownian path $f$.
\ep
}
\end{Remark}

We conclude this subsection with an integration by parts formula for the deterministic forward and backward integrals, whose simple proof is omitted.

\begin{Proposition}
\label{P:BVI}
Let $f\colon[a,b]\rightarrow\R$ be a c\`{a}dl\`{a}g function and $g\colon[a,b]\rightarrow\R$ be a bounded variation function. Then, the following \textbf{integration by parts formulae} hold:
\begin{align}
\int_{[a,b]}  g(s) d^-f(s) \ &= \ g(b) f(b)  - \int_{]a,b]} f(s) dg(s), \label{E:IbyP-} \\
\int_{[a,b]}  g(s) d^+f(s) \ &= \ g(b) f(b^-)  - \int_{]a,b]} f(s^-) dg(s). \label{E:IbyP}
\end{align}
\end{Proposition}

\subsection{The spaces $\mathscr C([-T,0])$ and $\mathscr C([-T,0[)$}
\label{SubS:Preliminaries}

Let $C([-T,0])$ denote the set of real continuous functions on $[-T,0]$, endowed with supremum norm $\|\eta\|_\infty = \sup_{x\in[-T,0]}|\eta(x)|$, for any $\eta\in C([-T,0])$.

\begin{Remark} 
{\rm
We shall develop functional It\^o calculus via regularization firstly for time-indepen\-dent functionals $\Uc\colon C([-T,0])\rightarrow\R$, since we aim at emphasizing that in our framework the time variable and the path play two distinct roles, as emphasized in the introduction. This, also, allows us to focus only on the definition of horizontal and vertical derivatives. Clearly, everything can be extended in an obvious way to the time-dependent case $\Uc\colon[0,T]\times C([-T,0])\rightarrow\R$, as we shall illustrate later.
\ep
}
\end{Remark}

Consider a map $\Uc\colon C([-T,0])\rightarrow\R$. Our aim is to derive a functional It\^o's formula for $\Uc$. To do this, we are led to define, in the spirit of \cite{dupire} and \cite{contfournie10}, the functional (i.e., horizontal and vertical) derivatives for $\Uc$. Since the definition of functional derivatives necessitates of discontinuous paths, in \cite{contfournie10} the idea is to consider functionals defined on the space of c\`adl\`ag trajectories $\D([-T,0])$. However, we can not, in general, extend in a unique way a functional $\Uc$ defined on $C([-T,0])$ to $\D([-T,0])$. Our idea, instead, is to consider a larger space than $C([-T,0])$, denoted by $\mathscr C([-T,0])$, which is the space of bounded trajectories on $[-T,0]$, continuous on $[-T,0[$ and with possibly a jump at $0$. We endow $\mathscr C([-T,0])$ with a (inductive) topology such that $C([-T,0])$ is dense in $\mathscr C([-T,0])$ with respect to this topology. Therefore, if $\Uc$ is continuous with respect to the topology of $\mathscr C([-T,0])$, then if it admits a continuous extension $u\colon \mathscr C([-T,0])\rightarrow\R$ this is necessarily unique.

\begin{Definition}
\label{D:scrC}
We denote by $\mathscr C([-T,0])$ the set of bounded functions $\eta\colon[-T,0]\rightarrow\R$ such that $\eta$ is continuous on $[-T,0[$, equipped with the topology we now describe.\\
\textbf{Convergence.} We endow $\mathscr C([-T,0])$ with a topology inducing the following convergence: $(\eta_n)_n$ converges to $\eta$ in $\mathscr C([-T,0])$ as $n$ tends to infinity if the following holds.
\begin{enumerate}
\item[\textup{(i)}] $\|\eta_n\|_\infty \leq C$, for any $n\in\N$, for some positive constant $C$ independent of $n$$;$
\item[\textup{(ii)}] $\sup_{x\in K}|\eta_n(x)-\eta(x)|\rightarrow0$ as $n$ tends to infinity, for any compact set $K\subset[-T,0[$$;$
\item[\textup{(iii)}] $\eta_n(0)\rightarrow\eta(0)$ as $n$ tends to infinity.
\end{enumerate}
\textbf{Topology.} For each compact $K\subset [-T,0[$ define the seminorm $p_K$ on $\mathscr C([-T,0])$ by
\[
p_K(\eta) \ = \ \sup_{x\in K}|\eta(x)| + |\eta(0)|, \qquad \forall\,\eta\in\mathscr C([-T,0]).
\]
Let $M>0$ and $\mathscr C_M([-T,0])$ be the set of functions in $\mathscr C([-T,0])$ which are bounded by $M$. Still denote $p_K$ the restriction of $p_K$ to $\mathscr C_M([-T,0])$ and consider the topology on $\mathscr C_M([-T,0])$ induced by the collection of seminorms $(p_K)_K$. Then, we endow $\mathscr C([-T,0])$ with the smallest topology $($inductive topology$)$ turning all the inclusions $i_M\colon \mathscr C_M([-T,0])\rightarrow\mathscr C([-T,0])$ into continuous maps.
\end{Definition}

\begin{Remark}
\label{R:Density}
{\rm
(i) Notice that $C([-T,0])$ is dense in $\mathscr C([-T,0])$, when endowed with the topology of $\mathscr C([-T,0])$. As a matter of fact, let $\eta\in\mathscr C([-T,0])$ and define, for any $n\in\N\backslash\{0\}$,
\[
\varphi_n(x)=
\begin{cases}
\eta(x), \qquad &-T\leq x\leq-1/n, \\
n(\eta(0)-\eta(-1/n))x + \eta(0), &-1/n<x\leq0.
\end{cases}
\]
Then, we see that $\varphi_n\in C([-T,0])$ and $\varphi_n\rightarrow\eta$ in $\mathscr C([-T,0])$.\\
Now, for any $a\in\R$ define
\begin{align*}
C_a([-T,0]) \ &:= \ \{\eta\in C([-T,0])\colon\eta(0)=a\}, \\\mathscr C_a([-T,0]) \ &:= \ \{\eta\in\mathscr C([-T,0])\colon\eta(0)=a\}.
\end{align*}
Then, $C_a([-T,0])$ is dense in $\mathscr C_a([-T,0])$ with respect to the topology of $\mathscr C([-T,0])$. \\
(ii) We provide two examples of functionals $\Uc\colon C([-T,0])\rightarrow\R$, continuous with respect to the topology of $\mathscr C([-T,0])$, and necessarily with respect to the topology of $C([-T,0])$; the proof is straightforward and 
not reported.
\begin{enumerate}
\item[(a)] $\Uc(\eta) = g(\eta(t_1),\ldots,\eta(t_n))$, for all $\eta\in C([-T,0])$, with $-T\leq t_1<\cdots<t_n\leq0$ and $g\colon\R^n\rightarrow\R$ continuous.
\item[(b)] $\Uc(\eta) = \int_{[-T,0]}\varphi(x)d^-\eta(x)$, for all $\eta\in C([-T,0])$, with $\varphi\colon[0,T]\rightarrow\R$ a c\`adl\`ag bounded variation function. Concerning this example, keep in mind that, using the integration by parts formula, $\Uc(\eta)$ admits the representation \eqref{E:IbyP-}.
\end{enumerate}
(iii) Consider the functional $\Uc(\eta) = 
\sup_{x\in[-T,0]}\eta(x)$, for all $\eta\in C([-T,0])$.
It is obviously continuous, but it
 is not continuous with 
respect to the topology of $\mathscr C([-T,0])$. As a matter of fact, for any $n\in\N$ consider $\eta_n\in C([-T,0])$ given by
\[
\eta_n(x) \ = \
\begin{cases}
0, \qquad\qquad\qquad &-T\leq x\leq-\frac{T}{2^n}, \\
\frac{2^{n+1}}{T}x+2, &-\frac{T}{2^n}<x\leq-\frac{T}{2^{n+1}}, \\
-\frac{2^{n+1}}{T}x, &-\frac{T}{2^{n+1}}<x\leq0.
\end{cases}
\]
Then, $\Uc(\eta_n)=\sup_{x\in[-T,0]}\eta_n(x)=1$, for any $n$. However, $\eta_n$ converges to the zero function in $\mathscr C([-T,0])$, as $n$ tends to infinity. This example plays an important role in the companion paper \cite{cosso_russo15b} to justify a weaker notion of solution to the path-dependent semilinear Kolmogorov equation.
\ep
}
\end{Remark}

To define the functional derivatives, we shall need to separate the
 ``past'' from the ``present'' of $\eta\in\mathscr C([-T,0])$. Indeed, roughly speaking, the horizontal derivative calls in the past values of $\eta$, namely $\{\eta(x)\colon x\in[-T,0[\}$, while the vertical derivative calls in the present value of $\eta$, namely $\eta(0)$. To this end, it is useful to introduce the space $\mathscr C([-T,0[)$.

\begin{Definition}
\label{D:scrC2}
We denote by $\mathscr C([-T,0[)$ the set of real-valued bounded continuous functions $\gamma\colon[-T,0[\rightarrow\R$, equipped with the topology we now describe.\\
\textbf{Convergence.} We endow $\mathscr C([-T,0[)$ with a topology inducing the following convergence: $(\gamma_n)_n$ converges to $\gamma$ in $\mathscr C([-T,0[)$ as $n$ tends to infinity if the following holds.
\begin{enumerate}
\item[\textup{(i)}] $\sup_{x\in[-T,0[}|\gamma_n(x)| \leq C$, for any $n\in\N$, for some positive constant $C$ independent of $n$$;$
\item[\textup{(ii)}] $\sup_{x\in K}|\gamma_n(x)-\gamma(x)|\rightarrow0$ as $n$ tends to infinity, for any compact set $K\subset[-T,0[$.
\end{enumerate}
\textbf{Topology.} For each compact $K\subset [-T,0[$ define the seminorm $q_K$ on $\mathscr C([-T,0[)$ by
\[
q_K(\gamma) \ = \ \sup_{x\in K}|\gamma(x)|, \qquad \forall\,\gamma\in\mathscr C([-T,0[).
\]
Let $M>0$ and $\mathscr C_M([-T,0[)$ be the set of functions in $\mathscr C([-T,0[)$ which are bounded by $M$. Still denote $q_K$ the restriction of $q_K$ to $\mathscr C_M([-T,0[)$ and consider the topology on $\mathscr C_M([-T,0[)$ induced by the collection of seminorms $(q_K)_K$. Then, we endow $\mathscr C([-T,0[)$ with the smallest topology $($inductive topology$)$ turning all the inclusions $i_M\colon \mathscr C_M([-T,0[)\rightarrow\mathscr C([-T,0[)$ into continuous maps.
\end{Definition}

\begin{Remark}
\label{R:Isomorphism}
{\rm
(i) Notice that $\mathscr C([-T,0])$ is isomorphic to $\mathscr C([-T,0[)\times\R$. As a matter of fact, it is enough to consider the map
\begin{align*}
J \colon \mathscr C([-T,0]) &\rightarrow \mathscr C([-T,0[)\times\R \\
\eta &\mapsto (\eta_{|[-T,0[},\eta(0)).
\end{align*}
Observe that $J^{-1}\colon \mathscr C([-T,0[)\times\R \rightarrow \mathscr C([-T,0])$ is given by $J^{-1}(\gamma,a) = \gamma1_{[-T,0[} + a1_{\{0\}}$.\\
(ii) $\mathscr C([-T,0])$ is a space which contains $C([-T,0])$ as a subset and it has the property of 
separating ``past'' from ``present''. Another space having the same property is $L^2([-T,0]; d \mu)$ where $\mu$ is the  sum of
the Dirac measure at zero and Lebesgue measure. Similarly as for item (i), that space is isomorphic to $L^2([-T,0]) \times \R $, which is a very popular space appearing in the analysis of functional dependent (as delay) equations, starting from \cite{cho}.
\ep
}
\end{Remark}

For every $u\colon \mathscr C([-T,0])\rightarrow\R$, we can now exploit the space $\mathscr C([-T,0[)$ to define a map $\tilde u\colon \mathscr C([-T,0[)\times\R\rightarrow\R$ where ``past'' and ``present'' are separated.

\begin{Definition}
\label{D:tildeu}
Let $u\colon \mathscr C([-T,0])\rightarrow\R$ and define $\tilde u\colon \mathscr C([-T,0[)\times\R\rightarrow\R$ as
\begin{equation}
\label{E:tildeu}
\tilde u(\gamma,a) \ := \ u(\gamma 1_{[-T,0[} + a1_{\{0\}}), \qquad \forall\,(\gamma,a)\in \mathscr C([-T,0[)\times\R.
\end{equation}
In particular, we have $u(\eta) = \tilde u(\eta_{|[-T,0[},\eta(0))$, for all $\eta\in \mathscr C([-T,0])$.
\end{Definition}

We conclude this subsection with a characterization of the dual spaces of $\mathscr C([-T,0])$ and $\mathscr C([-T,0[)$, which has an independent interest. Firstly, we need to introduce the set $\Mc([-T,0])$ of finite signed Borel measures on $[-T,0]$. We also denote $\Mc_0([-T,0])\subset\Mc([-T,0])$ the set of measures $\mu$ such that $\mu(\{0\})=0$.

\begin{Proposition}
\label{P:DualSpace}
Let $\Lambda\in \mathscr C([-T,0])^*$, the dual space of $\mathscr C([-T,0])$. Then, there exists a unique $\mu\in\Mc([-T,0])$ such that
\[
\Lambda\eta \ = \ \int_{[-T,0]} \eta(x) \mu(dx), \qquad \forall\,\eta\in \mathscr C([-T,0]).
\]
\end{Proposition}
\textbf{Proof.}
Let $\Lambda\in \mathscr C([-T,0])^*$ and define
\[
\tilde\Lambda\varphi \ := \ \Lambda\varphi, \qquad \forall\,\varphi\in C([-T,0]).
\]
Notice that $\tilde\Lambda\colon C([-T,0])\rightarrow\R$ is a continuous functional on the Banach space $C([-T,0])$ endowed with the supremum norm $\|\cdot\|_\infty$. Therefore $\tilde\Lambda\in C([-T,0])^*$ and it follows from Riesz representation theorem (see, e.g., Theorem 6.19 in \cite{rudin}) that there exists a unique $\mu\in\Mc([-T,0])$ such that
\[
\tilde\Lambda\varphi \ = \ \int_{[-T,0]} \varphi(x)\mu(dx), \qquad \forall\,\varphi\in C([-T,0]).
\]
Obviously $\tilde\Lambda$ is also continuous with respect to the topology of $\mathscr C([-T,0])$. Since $C([-T,0])$ is dense in $\mathscr C([-T,0])$ with respect to the topology of $\mathscr C([-T,0])$, we deduce that there exists a unique continuous extension of $\tilde\Lambda$ to $\mathscr C([-T,0])$, which is clearly given by
\[
\Lambda\eta \ = \ \int_{[-T,0]} \eta(x)\mu(dx), \qquad \forall\,\eta\in \mathscr C([-T,0]).
\]
\ep

\begin{Proposition}
\label{P:DualSpace[}
Let $\Lambda\in \mathscr C([-T,0[)^*$, the dual space of $\mathscr C([-T,0[)$. Then, there exists a unique $\mu\in\Mc_0([-T,0])$ such that
\[
\Lambda\gamma \ = \ \int_{[-T,0[} \gamma(x) \mu(dx), \qquad \forall\,\gamma\in \mathscr C([-T,0[).
\]
\end{Proposition}
\textbf{Proof.}
Let $\Lambda\in \mathscr C([-T,0[)^*$ and define
\begin{equation}
\label{E:DualSpace[Proof}
\tilde\Lambda\eta \ := \ \Lambda(\eta_{|[-T,0[}), \qquad \forall\,\eta\in\mathscr C([-T,0]).
\end{equation}
Notice that $\tilde\Lambda\colon\mathscr C([-T,0])\rightarrow\R$ is a continuous functional on $\mathscr C([-T,0])$. It follows from Proposition \ref {P:DualSpace} that there exists a unique $\mu\in\Mc([-T,0])$ such that
\begin{equation}
\label{E:DualSpace[Proof2}
\tilde\Lambda\eta \ = \ \int_{[-T,0]} \eta(x)\mu(dx) \ = \ \int_{[-T,0[} \eta(x)\mu(dx) + \eta(0)\mu(\{0\}), \qquad \forall\,\eta\in\mathscr C([-T,0]).
\end{equation}
Let $\eta_1,\eta_2\in\mathscr C([-T,0])$ be such that $\eta_1 1_{[-T,0[}=\eta_2 1_{[-T,0[}$. Then, we see from \eqref{E:DualSpace[Proof} that $\tilde\Lambda\eta_1=\tilde\Lambda\eta_2$, which in turn implies from \eqref{E:DualSpace[Proof2} that $\mu(\{0\})=0$. In conclusion, $\mu\in\Mc_0([-T,0])$ and $\Lambda$ is given by
\[
\Lambda\gamma \ = \ \int_{[-T,0[} \gamma(x)\mu(dx), \qquad \forall\,\gamma\in \mathscr C([-T,0[).
\]
\ep

\subsection{Functional derivatives and functional It\^o's formula}
\label{SubS:PathwiseDerivatives}

In the present section we shall prove one of the main result of this section, namely the functional It\^o's formula for $\Uc\colon C([-T,0])\rightarrow\R$ and, more generally, for $\Uc\colon[0,T]\times C([-T,0])\rightarrow\R$. We begin introducing the functional derivatives in the spirit of Dupire \cite{dupire}, firstly for a functional $u\colon\mathscr C([-T,0])\rightarrow\R$, and then for $\Uc\colon C([-T,0])\rightarrow\R$.

\begin{Definition}
\label{D:DupireDerivatives}
\quad\\ Consider $u\colon\mathscr C([-T,0])\rightarrow\R$ and $\eta\in\mathscr C([-T,0])$. \\
\textup{(i)} We say that $u$ admits \textbf{horizontal derivative} at $\eta$ if the following limit exists and it is finite:
\[
D^H u(\eta) \ := \ \lim_{\eps\rightarrow0^+} \frac{u(\eta(\cdot)1_{[-T,0[}+\eta(0)1_{\{0\}}) - u(\eta(\cdot-\eps)1_{[-T,0[}+\eta(0)1_{\{0\}})}{\eps}.
\]
\textup{(i)'} Let $\tilde u$ be as in \eqref{E:tildeu}. We say that $\tilde u$ admits \textbf{horizontal derivative} at $(\gamma,a)\in\mathscr C([-T,0[)\times\R$ if the following limit exists and it is finite:
\begin{equation}
\label{E:DHtildeu}
D^H\tilde u(\gamma,a) \ := \ \lim_{\eps\rightarrow0^+} \frac{\tilde u(\gamma(\cdot),a) - \tilde u(\gamma(\cdot-\eps),a)}{\eps}.
\end{equation}
Notice that if $D^H u(\eta)$ exists then $D^H\tilde u(\eta_{|[-T,0[},\eta(0))$ exists and they are equal; viceversa, whenever $D^H\tilde u(\gamma,a)$ exists then $D^H u(\gamma 1_{[-T,0[} + a 1_{\{0\}})$ exists and they are equal. \\
\textup{(ii)} We say that $u$ admits \textbf{first-order vertical
 derivative} at $\eta$ if the first-order partial derivative at $(\eta_{|[-T,0[},\eta(0))$ of $\tilde u$ with respect to its second argument, which we denote by $\partial_a\tilde u(\eta_{|[-T,0[},\eta(0))$, exists and we set
\[
D^V u(\eta) \ := \ \partial_a \tilde u(\eta_{|[-T,0[},\eta(0)).
\]
\textup{(iii)} We say that $u$ admits \textbf{second-order vertical derivative} at $\eta$ if the second-order partial derivative at $(\eta _{|[-T,0[},\eta(0))$ of $\tilde u$ with respect to its second argument, which we denote by $\partial_{aa}^2\tilde u(\eta_{|[-T,0[},\eta(0))$, exists and we set
\[
D^{VV} u(\eta) \ := \ \partial_{aa}^2 \tilde u(\eta_{|[-T,0[},\eta(0)).
\]
\end{Definition}

\begin{Definition}
\label{D:C12}
We say that $u\colon \mathscr C([-T,0])\rightarrow\R$ is of class $\mathscr C^{1,2}(\textup{past}\times\textup{present})$ if the following holds.
\begin{enumerate}
\item[\textup{(i)}] $u$ is continuous$;$
\item[\textup{(ii)}] $D^H u$ exists everywhere on $\mathscr C([-T,0])$ and for every $\gamma\in \mathscr C([-T,0[)$ the map
\[
(\eps,a)\longmapsto D^H\tilde u(\gamma(\cdot-\eps),a), \qquad (\eps,a)\in[0,\infty[\times\R
\]
is continuous on $[0,\infty[\times\R$$;$
\item[\textup{(iii)}] $D^V u$ and $D^{VV}u$ exist everywhere on 
$\mathscr C([-T,0])$ and are continuous.
\end{enumerate}
\end{Definition}

\begin{Remark}
{\rm
Notice that in Definition \ref{D:C12} we still obtain the same class of functions $\mathscr C^{1,2}(\textup{past}\times\textup{present})$ if we substitute point (ii) with:
\begin{enumerate}
\item[\textup{(ii')}] $D^H u$ exists everywhere on $\mathscr C([-T,0])$ and for every $\gamma\in \mathscr C([-T,0[)$ there exists 
$\delta(\gamma) \in ]0,\infty]$ such that the map
\begin{equation}
\label{E:map}
(\eps,a)\longmapsto D^H\tilde u(\gamma(\cdot-\eps),a), \qquad (\eps,a)\in[0,\infty[\times\R
\end{equation}
is continuous on $[0,\delta(\gamma)[\times\R$.
\end{enumerate}
In particular, if (ii') holds then we can always take $\delta(\gamma) = \infty$ for any $\gamma\in \mathscr C([-T,0[)$, which implies (ii). To prove this last statement, let us proceed by contradiction assuming that
\[
\delta^*(\gamma) \ = \ \sup\big\{\delta(\gamma)>0\colon\text{the map \eqref{E:map} is continuous on }[0,\delta(\gamma)[\times\R\big\} \ < \ \infty.
\]
Notice that $\delta^*(\gamma)$ is in fact a \emph{max}, therefore the map \eqref{E:map} is continuous on $[0,\delta^*(\gamma)[\times\R$. Now, define $\bar{\gamma}(\cdot) := \gamma(\cdot-\delta^*(\gamma))$. Then, by condition (ii') there exists $\delta(\bar{\gamma})>0$ such that the map
\[
(\eps,a)\longmapsto D^H\tilde u(\bar{\gamma}(\cdot-\eps),a) = D^H\tilde u(\gamma(\cdot-\eps-\delta^*(\gamma)),a)
\]
is continuous on $[0,\delta(\bar{\gamma})[\times\R$. This shows that the map \eqref{E:map} is continuous on $[0,\delta^*(\gamma)+\delta(\bar{\gamma})[\times\R$, a contradiction with the definition of $\delta^*(\gamma)$.
\ep
}
\end{Remark}

\noindent We can now provide the definition of functional derivatives for a map $\Uc\colon$ $C([-T,0])$ $\rightarrow$ $\R$.

\begin{Definition}
\label{D:DupireDerivativesUc}
Let $\Uc\colon C([-T,0])\rightarrow\R$ and $\eta\in C([-T,0])$. Suppose that there exists a unique extension $u\colon\mathscr C([-T,0])\rightarrow\R$ of $\Uc$ $($e.g., if $\Uc$ is continuous with respect to the topology of $\mathscr C([-T,0])$$)$. Then we define:\\
\textup{(i)} The \textbf{horizontal derivative} of $\Uc$ at $\eta$ as:
\[
D^H \Uc(\eta) \ := \ D^H u(\eta).
\]
\textup{(ii)} The \textbf{first-order vertical derivative} of $\Uc$ at $\eta$ as:
\[
D^V \Uc(\eta) \ := \ D^V u(\eta).
\]
\textup{(iii)} The \textbf{second-order vertical derivative} of $\Uc$ at $\eta$ as:
\[
D^{VV} \Uc(\eta) \ := \ D^{VV} u(\eta).
\]
\end{Definition}

\begin{Definition}
\label{D:C12Uc}
We say that $\Uc\colon C([-T,0])\rightarrow\R$ is $C^{1,2}(\textup{past}\times\textup{present})$ if $\Uc$ admits a $($necessarily unique$)$ extension $u\colon \mathscr C([-T,0])\rightarrow\R$ of class $\mathscr C^{1,2}(\textup{past}\times\textup{present})$.
\end{Definition}

\begin{Theorem}
\label{T:Ito}
Let $\Uc\colon C([-T,0])\rightarrow\R$ be of class $C^{1,2}(\textup{past}\times\textup{present})$ and $X=(X_t)_{t\in[0,T]}$ be a real 
continuous finite quadratic variation process. Then, the following \textbf{functional It\^o formula} holds, $\P$-a.s.,
\begin{equation}
\label{E:Ito}
\Uc(\X_t) \ = \ \Uc(\X_0) + \int_0^t D^H \Uc(\X_s)ds + \int_0^t D^V \Uc(\X_s) d^- X_s + \frac{1}{2}\int_0^t D^{VV}\Uc(\X_s)d[X]_s,
\end{equation}
for all $0 \leq t \leq T$, where the window process $\X$ was defined in \eqref{WindowProcess}.
\end{Theorem}
\textbf{Proof.}
Fix $t\in[0,T]$ and consider the quantity
\[
I_0(\eps,t) \ = \ \int_0^t \frac{\Uc(\X_{s+\eps}) - \Uc(\X_s)}{\eps} ds \ = \ \frac{1}{\eps} \int_t^{t+\eps} \Uc(\X_s) ds - \frac{1}{\eps} \int_0^\eps \Uc(\X_s) ds, \qquad \eps>0.
\]
Since $(\Uc(\X_s))_{s\geq0}$ is continuous, we have that $I_0(\eps,t)$ converges ucp to $\Uc(\X_t) - \Uc(\X_0)$, namely $\sup_{0 \leq t \leq T}|I_0(\eps,t)-(\Uc(\X_t) - \Uc(\X_0))|$ converges to zero in probability when $\eps\rightarrow0^+$. On the other hand, we can write $I_0(\eps,t)$ in terms of the function $\tilde u$, defined in \eqref{E:tildeu}, as follows:
\[
I_0(\eps,t) \ = \ \int_0^t \frac{\tilde u(\X_{s+\eps|[-T,0[},X_{s+\eps}) - \tilde u(\X_{s|[-T,0[},X_s)}{\eps} ds.
\]
Now we split $I_0(\eps,t)$ into two terms:
\begin{align}
I_1(\eps,t) \ &= \ \int_0^t \frac{\tilde u(\X_{s+\eps|[-T,0[},X_{s+\eps}) - \tilde u(\X_{s|[-T,0[},X_{s+\eps})}{\eps} ds, \label{E:I1} \\
I_2(\eps,t) \ &= \ \int_0^t \frac{\tilde u(\X_{s|[-T,0[},X_{s+\eps}) - \tilde u(\X_{s|[-T,0[},X_s)}{\eps} ds. \label{E:I2}
\end{align}
We begin proving that
\begin{equation}
\label{E:I1ucp}
I_1(\eps,t) \underset{\eps\rightarrow0^+}{\overset{\text{ucp}}{\longrightarrow}} \int_0^t D^H \Uc(\X_s) ds.
\end{equation}
Firstly, fix $\gamma\in \mathscr C([-T,0[)$ and define
\[
\phi(\eps,a) \ := \ \tilde u(\gamma(\cdot-\eps),a), \qquad (\eps,a)\in[0,\infty[\times\R.
\]
Then, denoting by $\partial_\eps^+ \phi$ the right partial derivative of $\phi$ with respect to $\eps$ and using formula \eqref{E:DHtildeu}, we find
\begin{align*}
\partial_\eps^+ \phi(\eps,a) \ &= \ \lim_{r\rightarrow0^+} \frac{\phi(\eps+r,a) - \phi(\eps,a)}{r} \\
&= \ -\lim_{r\rightarrow0^+} \frac{\tilde u(\gamma(\cdot-\eps),a) - \tilde u(\gamma(\cdot-\eps-r),a)}{r} \\
&= \ -D^H \tilde u(\gamma(\cdot-\eps),a), \qquad \forall\,(\eps,a)\in[0,\infty[\times\R.
\end{align*}
Since $u\in\mathscr C^{1,2}(\textup{past}\times\textup{present})$, we see from Definition \ref{D:C12}(ii), that $\partial_\eps^+ \phi$ is continuous on $[0,\infty[\times\R$. It follows from a standard differential calculus' result (see for example Corollary 1.2, Chapter 2, in \cite{pazy83}) that $\phi$ is continuously differentiable on $[0,\infty[\times\R$ with respect to its first argument. Then, for every $(\eps,a)\in[0,\infty[\times\R$, from the fundamental theorem of calculus, we have
\[
\phi(\eps,a) - \phi(0,a) \ = \ \int_0^\eps \partial_\eps\phi (r,a) dr,
\]
which in terms of $\tilde u$ reads
\begin{equation}
\label{E:tildeuFundThmCalc}
\tilde u(\gamma(\cdot),a) - \tilde u(\gamma(\cdot-\eps),a) \ = \ \int_0^\eps D^H \tilde u(\gamma(\cdot-r),a) dr.
\end{equation}
Now, we rewrite, by means of a shift in time, the term $I_1(\eps,t)$ in \eqref{E:I1} as follows:
\begin{align}
\label{E:I1bis}
I_1(\eps,t) \ &= \int_0^t \frac{\tilde u(\X_{s|[-T,0[},X_s) - \tilde u(\X_{s-\eps|[-T,0[},X_s)}{\eps} ds + \int_t^{t+\eps} \frac{\tilde u(\X_{s|[-T,0[},X_s) - \tilde u(\X_{s-\eps|[-T,0[},X_s)}{\eps} ds \notag \\
&\quad - \int_0^{\eps} \frac{\tilde u(\X_{s|[-T,0[},X_s) - \tilde u(\X_{s-\eps|[-T,0[},X_s)}{\eps} ds.
\end{align}
Plugging \eqref{E:tildeuFundThmCalc} into \eqref{E:I1bis}, setting $\gamma =  \X_s, a = X_s$, we obtain
\begin{align}
\label{E:I1R1}
I_1(\eps,t) \ &= \int_0^t \frac{1}{\eps} \bigg(\int_0^\eps D^H \tilde u(\X_{s-r|[-T,0[},X_s)dr\bigg)ds + \int_t^{t+\eps} \frac{1}{\eps} \bigg(\int_0^\eps D^H \tilde u(\X_{s-r|[-T,0[},X_s)dr\bigg)ds \notag \\
&\quad - \int_0^\eps \frac{1}{\eps} \bigg(\int_0^\eps D^H \tilde u(\X_{s-r|[-T,0[},X_s)dr\bigg)ds.
\end{align}
Observe that
\[
\int_0^t \frac{1}{\eps} \bigg(\int_0^\eps D^H \tilde u(\X_{s-r|[-T,0[},X_s)dr\bigg)ds \underset{\eps\rightarrow0^+}{\overset{\text{ucp}}{\longrightarrow}} \int_0^t D^H u(\X_s) ds.
\]
Similarly, we see that the other two terms in \eqref{E:I1R1} converge ucp to zero. As a consequence, we get \eqref{E:I1ucp}.

Regarding $I_2(\eps,t)$ in \eqref{E:I2}, it can be written, by means of the following standard Taylor's expansion for a function $f\in C^2(\R)$:
\begin{align*}
f(b) \ &= \ f(a)+f'(a)(b-a)+\frac{1}{2}f''(a)(b-a)^2 + \int_0^1(1-\alpha)\big(f''(a+\alpha(b-a))-f''(a)\big)(b-a)^2d\alpha,
\end{align*}
as the sum of the following three terms:
\begin{align*}
I_{21}(\eps,t) \ &= \ \int_0^t \partial_a \tilde u(\X_{s|[-T,0[},X_s) \frac{X_{s+\eps}-X_s}{\eps} ds \\
I_{22}(\eps,t) \ &= \ \frac{1}{2}\int_0^t \partial_{aa}^2 \tilde u(\X_{s|[-T,0[},X_s) \frac{(X_{s+\eps}-X_s)^2}{\eps} ds \\
I_{23}(\eps,t) \ &= \ \int_0^t \bigg(\int_0^1(1-\alpha) \big( \partial_{aa}^2 \tilde u(\X_{s|[-T,0[},X_s + \alpha(X_{s+\eps}-X_s)) \\
&\quad\;\, - \partial_{aa}^2 \tilde u(\X_{s|[-T,0[},X_s) \big) \frac{(X_{s+\eps}-X_s)^2}{\eps} d\alpha \bigg) ds.
\end{align*}
By similar arguments as in Proposition 1.2 of \cite{russovallois95}, we have 
\[
I_{22}(\eps,t) \underset{\eps\rightarrow0^+}{\overset{\text{ucp}}{\longrightarrow}} \frac{1}{2}\int_0^t\partial_{aa}^2 \tilde u(\X_{s|[-T,0[},X_s) d[X]_s = \frac{1}{2}\int_0^t D^{VV} u(\X_s) d[X]_s.
\]
Regarding $I_{23}(\eps,t)$, for every $\omega\in\Omega$, define $\psi_\omega\colon[0,T]\times[0,1]\times[0,1]\rightarrow\R$ as
\[
\psi_\omega(s,\alpha,\eps) \ := \ (1-\alpha) \partial_{aa}^2 \tilde u\big(\X_{s|[-T,0[}(\omega),X_s(\omega) + \alpha(X_{s+\eps}(\omega)-X_s(\omega))\big),
\]
for all $(s,\alpha,\eps)\in[0,T]\times[0,1]\times[0,1]$. Notice that $\psi_\omega$ is uniformly continuous. Denote $\rho_{\psi_\omega}$ its continuity modulus, then
\[
\sup_{t\in[0,T]} |I_{23}(\eps,t)| \ \leq \ \int_0^T\rho_{\psi_\omega}(\eps) \frac{(X_{s+\eps}-X_s)^2}{\eps}ds.
\]
Since $X$ has finite quadratic variation, we deduce that $I_{23}(\eps,t)\rightarrow0$ ucp as $\eps\rightarrow0^+$. Finally, because of $I_0(\eps,t)$, $I_1(\eps,t)$, $I_{22}(\eps,t)$, and $I_{23}(\eps,t)$ converge ucp, it follows that the forward integral exists:
\[
I_{21}(\eps,t) \underset{\eps\rightarrow0^+}{\overset{\text{ucp}}{\longrightarrow}} \int_0^t \partial_a \tilde u(\X_{s|[-T,0[},X_s) d^- X_s = \int_0^t D^V u(\X_s) d^- X_s,
\]
from which the claim follows.
\ep

\begin{Remark}
{\rm
We notice that, under the assumptions of Theorem \ref{T:Ito}, the forward integral $\int_0^t D^V\Uc(\X_s)d^-X_s$ exists as a ucp limit, which is generally not required.
\ep
}
\end{Remark}

We conclude this subsection providing the functional It\^o formula for a map $\Uc\colon[0,T]\times C([-T,0])$ $\rightarrow\R$ depending also on the time variable. Firstly, we notice that for a map $\Uc\colon[0,T]\times C([-T,0])\rightarrow\R$ (resp. $u\colon[0,T]\times\mathscr C([-T,0])\rightarrow\R$) the functional derivatives $D^H\Uc$, $D^V\Uc$, and $D^{VV}\Uc$ (resp. $D^Hu$, $D^Vu$, and $D^{VV}u$) are defined in an obvious way as in Definition \ref{D:DupireDerivativesUc} (resp. Definition \ref{D:DupireDerivatives}). Moreover, given $u\colon[0,T]\times\mathscr C([-T,0])\rightarrow\R$ we can define, as in Definition \ref{D:tildeu}, a map $\tilde u\colon[0,T]\times\mathscr C([-T,0[)\times\R\rightarrow\R$. Then, we can give the following definitions.

\begin{Definition}
\label{D:C12Time}
Let $I$ be $[0,T[$ or $[0,T]$. We say that $u\colon I\times\mathscr C([-T,0])\rightarrow\R$ is of class $\mathscr C^{1,2}((I\times\textup{past})\times\textup{present})$ if the properties below hold.
\begin{enumerate}
\item[\textup{(i)}] $u$ is continuous$;$
\item[\textup{(ii)}] $\partial_tu$ exists everywhere on $I\times\mathscr C([-T,0])$ and is continuous$;$
\item[\textup{(iii)}] $D^H u$ exists everywhere on $I\times\mathscr C([-T,0])$ and for every $\gamma\in \mathscr C([-T,0[)$ the map
\[
(t,\eps,a)\longmapsto D^H\tilde u(t,\gamma(\cdot-\eps),a), \qquad (t,\eps,a)\in I\times[0,\infty[\times\R
\]
is continuous on $I\times[0,\infty[\times\R$$;$
\item[\textup{(iv)}] $D^V u$ and $D^{VV}u$ exist everywhere on $I\times\mathscr C([-T,0])$ and are continuous.
\end{enumerate}
\end{Definition}

\begin{Definition}
\label{D:C12UcTime}
Let $I$ be $[0,T[$ or $[0,T]$. We say that $\Uc\colon I\times C([-T,0])\rightarrow\R$ is $C^{1,2}((I\times\textup{past})\times\textup{present}))$ if $\Uc$ admits a $($necessarily unique$)$ extension $u\colon I\times\mathscr C([-T,0])\rightarrow\R$ of class $\mathscr C^{1,2}((I\times\textup{past})\times\textup{present})$.
\end{Definition}

We can now state the functional It\^o formula, whose proof is not reported, since it can be done along the same lines as Theorem \ref{T:Ito}.

\begin{Theorem}
\label{T:ItoTime}
Let $\Uc\colon[0,T]\times C([-T,0])\rightarrow\R$ be of class $C^{1,2}(([0,T]\times\textup{past})\times\textup{present})$ and $X=(X_t)_{t\in[0,T]}$ be a real continuous finite quadratic variation process. Then, the following \textbf{functional It\^o formula} holds, $\P$-a.s.,
\begin{align}
\label{E:ItoTime}
\Uc(t,\X_t) \ &= \ \Uc(0,\X_0) + \int_0^t \big(\partial_t\Uc(s,\X_s) + D^H \Uc(s,\X_s)\big)ds + \int_0^t D^V \Uc(s,\X_s) d^- X_s \notag \\
&\quad \ + \frac{1}{2}\int_0^t D^{VV}\Uc(s,\X_s)d[X]_s,
\end{align}
for all $0 \leq t \leq T$.
\end{Theorem}

\begin{Remark}
{\rm
Notice that, as a particular case, choosing $\Uc(t,\eta)=F(t,\eta(0))$, for any $(t,\eta)\in[0,T]\times C([-T,0])$, with $F\in C^{1,2}([0,T]\times\R)$, we retrieve the classical It\^o's formula 
for finite quadratic variation processes, i.e. \eqref{E:ITOFQV}. More precisely, in this case $\Uc$ admits as unique continuous extension the map $u\colon[0,T]\times\mathscr C([-T,0])\rightarrow\R$ given by $u(t,\eta)=F(t,\eta(0))$, for all $(t,\eta)\in[0,T]\times\mathscr C([-T,0])$. Moreover, we see that $D^H\Uc\equiv0$, while $D^V\Uc=\partial_x F$ and $D^{VV}\Uc=\partial_{xx}^2 F$, where $\partial_x F$ (resp. $\partial_{xx}^2F$) denotes the first-order (resp. second-order) partial derivative of $F$ with respect to its second argument.
\ep
}
\end{Remark}

\subsection{Comparison with Banach space valued calculus via regularization}
\label{SubS:Comparison}

In the present subsection our aim is to make a link between functional It\^o calculus, as derived in this paper, and Banach space valued stochastic calculus via regularization for window processes, which has been conceived in \cite{DGR}, see also \cite{DGR2, digirrusso12, DGRnote},
and  \cite{digirfabbrirusso13} for more recent developments. More precisely, our purpose is to identify the building blocks of our functional It\^o's formula \eqref{E:Ito} with the terms appearing in the It\^o's formula derived in Theorem 6.3 and Section 7.2 in \cite{digirfabbrirusso13}. While it is expected that the vertical derivative $D^V\Uc$ can be identified with the term $D_{dx}^{\delta_0}\Uc$ of the Fr\'echet derivative, it is more difficult to guess to which terms the horizontal derivative $D^H\Uc$ corresponds. To clarify this latter point, in this subsection we derive two formulae which express $D^H\Uc$ in terms of Fr\'echet derivatives of $\Uc$.

\vspace{3mm}

Let us introduce some useful notations. We denote by $BV([-T,0])$ the set of c\`adl\`ag bounded variation functions on $[-T,0]$, which is a Banach space when equipped with the norm
\[
\|\eta\|_{BV([-T,0])} \ := \ |\eta(0)| + \|\eta\|_{\textup{Var}([-T,0])}, \qquad \eta\in BV([-T,0]),
\]
where $\|\eta\|_{\textup{Var}([-T,0])}=|d\eta|([-T,0])$ and $|d\eta|$ is the total variation measure associated to the measure $d\eta\in\Mc([-T,0])$ generated by $\eta$: $d\eta([-T,-t])=\eta(-t)-\eta(-T)$, $t\in[-T,0]$. We recall from subsection \ref{SubS:Background} that we extend $\eta\in BV([-T,0])$ to all $x\in\R$ setting $\eta(x) = 0$, $x<-T$, and $\eta(x)=\eta(0)$, $x\geq0$. Let us now introduce some useful facts about tensor products of Banach spaces.

\begin{Definition}
\label{D:Tensor}
Let $(E,\|\cdot\|_E)$ and $(F,\|\cdot\|_F)$ be two Banach spaces.\\
\textup{(i)} We shall denote by $E\otimes F$ the \textbf{algebraic tensor product} of $E$ and $F$, defined as the set of elements of the form $v = \sum_{i=1}^n e_i\otimes f_i$, for some positive integer $n$, where $e\in E$ and $f\in F$. The map $\otimes\colon E\times F\rightarrow E\otimes F$ is bilinear.\\
\textup{(ii)} We endow $E\otimes F$ with the \textbf{projective norm} $\pi$
defined as follows:
\[
\pi(v) \ := \ \inf\bigg\{\sum_{i=1}^n \|e_i\|_E\|f_i\|_F \ \colon \ v = \sum_{i=1}^n e_i\otimes f_i\bigg\}, \qquad \forall\,v\in E\otimes F.
\]
\textup{(iii)} We denote by $E\hat\otimes_\pi F$ the Banach space obtained as the completion of $E\otimes F$ for the norm $\pi$. We shall refer to $E\hat\otimes_\pi F$ as the \textbf{tensor product of the Banach spaces $E$ and $F$}.\\
\textup{(iv)} If $E$ and $F$ are Hilbert spaces, we denote $E\hat\otimes_h F$ the \textbf{Hilbert tensor product}, which is still a Hilbert space obtained as the completion of $E\otimes F$ for the scalar product $\langle e'\otimes f',e''\otimes f''\rangle := \langle e',e''\rangle_E\langle f',f''\rangle_F$, for any $e',e''\in E$ and $f',f''\in F$.\\
\textup{(v)} The symbols $E\hat\otimes_\pi^2$ and $e\otimes^2$ denote, respectively, the Banach space $E\hat\otimes_\pi E$ and the element $e\otimes e$ of the algebraic tensor product $E\otimes E$.
\end{Definition}

\begin{Remark}
\label{R:Identification}
{\rm
(i) The projective norm $\pi$ belongs to the class of the so-called \emph{reasonable crossnorms} $\alpha$ on $E\otimes F$, verifying $\alpha(e\otimes f)=\|e\|_E\|f\|_F$.\\
(ii) We notice, proceeding for example as in \cite{digirrusso12} (see, in particular, formula (2.1) in \cite{digirrusso12}; for more information on this subject we refer to \cite{ryan02}), that the dual $(E\hat\otimes_\pi F)^*$ of $E\hat\otimes_\pi F$ is isomorphic to the space of continuous bilinear forms $\Bc i(E,F)$, equipped with the norm $\|\cdot\|_{E,F}$ defined as
\[
\|\Phi\|_{E,F} \ := \ \sup_{\substack{e\in E, f\in F \\ \|e\|_E,\|f\|_F \leq 1}} |\Phi(e,f)|, \qquad \forall\,\Phi\in\Bc i(E,F).
\]
\ep
}
\end{Remark}

\begin{Definition}
\label{D:C12Frechet}
Let $E$ be a Banach space. We say that $\Uc\colon E\rightarrow\R$ is of class $C^2(E)$ if
\begin{enumerate}
\item[(i)] $D \Uc$, the first Fr\'echet derivative of $\Uc$, belongs to $C(E; E^*)$ and
\item[(ii)] $D^2 \Uc$, the second Fr\'echet derivative of $\Uc$, belongs to $C(E; \Bc i(E,E))$.
\end{enumerate}
\end{Definition}

\begin{Remark}
{\rm
Take $E = C([-T,0])$ in Definition \ref{D:C12Frechet}.\\
(i) \emph{First Fr\'echet derivative $D\Uc$.} We have
\[
D\Uc \colon C([-T,0]) \ \longrightarrow \ (C([-T,0]))^* \cong \Mc([-T,0]).
\]
For every $\eta\in C([-T,0])$, we shall denote $D_{dx}\Uc(\eta)$ the unique measure in $\Mc([-T,0])$ such that
\[
D\Uc(\eta)\varphi \ = \ \int_{[-T,0]} \varphi(x) D_{dx}\Uc(\eta), \qquad \forall\,\varphi\in C([-T,0]).
\]
Notice that $\Mc([-T,0])$ can be represented as the direct sum $\Mc([-T,0]) = \Mc_0([-T,0])\oplus\Dc_0$, where we recall that $\Mc_0([-T,0])$ is the subset of $\Mc([-T,0])$ of measures $\mu$ such that $\mu(\{0\})=0$, instead $\Dc_0$ (which is a shorthand for $\Dc_0([-T,0])$) denotes the one-dimensional space of measures which are multiples of the Dirac measure $\delta_0$. For every $\eta\in C([-T,0])$ we denote by $(D_{dx}^\perp \Uc(\eta),D_{dx}^{\delta_0}\Uc(\eta))$ the unique pair in $\Mc_0([-T,0])\oplus\Dc_0$ such that
\[
D_{dx}\Uc(\eta) \ = \ D_{dx}^\perp \Uc(\eta) + D_{dx}^{\delta_0}\Uc(\eta).
\]
(ii) \emph{Second Fr\'echet derivative $D^2\Uc$.} We have
\begin{align*}
D^2\Uc \colon C([-T,0]) \ \longrightarrow 
\Bc i(C([-T,0]),C([-T,0])) 
\cong (C([-T,0])\hat\otimes_\pi C([-T,0]))^*,
\end{align*}
where we used the identifications of Remark \ref{R:Identification}(iii). Consider $\eta\in C([-T,0])$; then a typical situation arises when there exists $D_{dx\,dy} \Uc(\eta)$ in $\Mc([-T,0]^2)$ for which $D^2\Uc(\eta)\in \Bc i(C([-T,0]),$ $C([-T,0]))$ admits the representation
\[
D^2\Uc(\eta)(\varphi,\psi) \ = \ \int_{[-T,0]^2} \varphi(x)\psi(y) D_{dx\,dy}\Uc(\eta), \qquad \forall\,\varphi,\psi\in C([-T,0]).
\]
Moreover, $D_{dx\,dy}\Uc(\eta)$ is uniquely determined.
\ep
}
\end{Remark}
The definition below was given in \cite{DGR}.
\begin{Definition}
\label{D:Chi-subspace}
Let $E$ be a Banach space. A Banach subspace $(\chi,\|\cdot\|_\chi)$ continuously injected into $(E\hat\otimes_\pi^2)^*$, i.e., $\|\cdot\|_\chi\geq\|\cdot\|_{(E\hat\otimes_\pi^2)^*}$, will be called a \textbf{Chi-subspace} $($of $(E\hat\otimes_\pi^2)^*$$)$.
\end{Definition}

\begin{Remark}
\label{R:Chi-subspace}
{\rm
Take $E=C([-T,0])$ in Definition \ref{D:Chi-subspace}. 
As indicated in \cite{DGR}, a typical example of Chi-subspace of $C([-T,0])\hat\otimes_\pi^2$ is $\Mc([-T,0]^2)$ equipped with the usual total variation norm, denoted by $\|\cdot\|_{\text{Var}}$. Another important Chi-subspace of $C([-T,0])\hat\otimes_\pi^2$ is the following, which is also a Chi-subspace of $\Mc([-T,0]^2)$:
\begin{align*}
\chi_0 \ &:= \ \big\{\mu\in\Mc([-T,0]^2)\colon\mu(dx,dy) = g_1(x,y)dxdy + \lambda_1\delta_0(dx)\otimes\delta_0(dy) \\
&\quad \ \ \, + g_2(x)dx\otimes\lambda_2\delta_0(dy) + \lambda_3\delta_0(dx)\otimes g_3(y)dy + g_4(x)\delta_y(dx)\otimes dy, \\
&\quad \ \ \; \, g_1\in L^2([-T,0]^2),\,g_2,g_3\in L^2([-T,0]),\,g_4\in L^\infty([-T,0]),\,\lambda_1,\lambda_2,\lambda_3\in\R\big\}.
\end{align*}
Using the notations of Example 3.4 and Remark 3.5 in \cite{digirrusso12}, to which we refer for more details on this subject, we notice that $\chi_0$ is indeed given by the direct sum $\chi_0 = L^2([-T,0]^2) \oplus \big(L^2([-T,0])\hat\otimes_h \Dc_0\big) \oplus \big(\Dc_0 \hat\otimes_h L^2([-T,0])\big) \oplus \Dc_{0,0}([-T,0]^2) \oplus Diag([-T,0]^2)$. In the sequel, we shall refer to the term $g_4(x)\delta_y(dx)\otimes dy$ as the {\bf diagonal component} and to $g_4(x)$ as the {\bf diagonal element} of $\mu$. 
\ep
}
\end{Remark}

\noindent We can now state our first representation result for $D^H\Uc$.

\begin{Proposition}
\label{P:DH=Dacdeta}
Let $\Uc\colon C([-T,0])\rightarrow\R$ be continuously Fr\'{e}chet differentiable. Suppose the following.
\begin{enumerate}
\item[\textup{(i)}] For any $\eta\in C([-T,0])$ there exists $D_x^{\textup{ac}}\Uc(\eta)\in BV([-T,0])$ such that
\[
D_{dx}^\perp \Uc(\eta) \ = \ D_x^{\textup{ac}}\Uc(\eta)dx.
\]
\item[\textup{(ii)}] There exist continuous extensions 
$($necessarily unique$)$
\[
u\colon\mathscr C([-T,0])\rightarrow\R, \qquad\qquad D_x^{\textup{ac}}u\colon\mathscr C([-T,0])\rightarrow BV([-T,0])
\]
of $\Uc$ and $D_x^{\textup{ac}}\Uc$, respectively.
\end{enumerate}
Then, for any $\eta\in C([-T,0])$,
\begin{equation}
\label{E:DH=Dacdeta}
D^H \Uc(\eta) \ = \ \int_{[-T,0]} D_x^{\textup{ac}} \Uc(\eta) d^+ \eta(x),
\end{equation}
where we recall that previous deterministic integral has been defined
in Section \ref{S211}.
In particular, the horizontal derivative $D^H \Uc(\eta)$ and the backward integral in \eqref{E:DH=Dacdeta} exist.
\end{Proposition}

\noindent\textbf{Proof.}
Let $\eta\in C([-T,0])$, then starting from the left-hand side of \eqref{E:DH=Dacdeta}, using the definition of $D^H\Uc(\eta)$, we are led to consider the following increment for the function $u$:
\begin{equation}
\label{E:FirstOrderDH}
\frac{u(\eta) - u(\eta(\cdot-\eps)1_{[-T,0[}+\eta(0)1_{\{0\}})}{\eps}.
\end{equation}
We shall expand \eqref{E:FirstOrderDH} using a Taylor formula. Firstly, notice that, since $\Uc$ is $C^1$ Fr\'echet on $C([-T,0])$, for every $\eta_1\in C([-T,0])$, with $\eta_1(0)=\eta(0)$,  from the fundamental theorem of calculus we have
\[
\Uc(\eta) - \Uc(\eta_1) \ = \ \int_0^1 \bigg( \int_{-T}^0 D_x^{\text{ac}}\Uc(\eta + \lambda(\eta_1-\eta))(\eta(x)-\eta_1(x)) dx\bigg) d\lambda.
\]
Recalling from Remark \ref{R:Density} the density of $C_{\eta(0)}([-T,0])$ in $\mathscr C_{\eta(0)}([-T,0])$ with respect to the topology of $\mathscr C([-T,0])$, we deduce the following Taylor's formula for $u$:
\begin{equation}
\label{E:u-u1}
u(\eta) - u(\eta_1) \ = \ \int_0^1 \bigg( \int_{-T}^0 D_x^{\text{ac}}u(\eta + \lambda(\eta_1-\eta))(\eta(x)-\eta_1(x)) dx\bigg) d\lambda,
\end{equation}
for all $\eta_1\in\mathscr C_{\eta(0)}([-T,0])$. As a matter of fact, for any $\delta\in]0,T/2]$ let (similarly to Remark \ref{R:Density}(i))
\[
\eta_{1,\delta}(x) \ := \
\begin{cases}
\eta_1(x), \qquad &-T\leq x\leq-\delta, \\
\frac{1}{\delta}(\eta_1(0)-\eta_1(-\delta))x + \eta_1(0), &-\delta<x\leq0
\end{cases}
\]
and $\eta_{1,0}:=\eta_1$. Then $\eta_{1,\delta}\in C([-T,0])$, for any $\delta\in]0,T/2]$, and $\eta_{1,\delta}\rightarrow\eta_1$ in $\mathscr C([-T,0])$, as $\delta\rightarrow0^+$. Now, define $f\colon[-T,0]\times[0,1]\times[0,T/2]\rightarrow\R$ as follows
\[
f(x,\lambda,\delta) \ := \ D_x^{\text{ac}}u(\eta + \lambda(\eta_{1,\delta}-\eta))(\eta(x)-\eta_{1,\delta}(x)),
\]
for all $(x,\lambda,\delta)\in[-T,0]\times[0,1]\times[0,T/2]$. Notice that $f$ is continuous and hence bounded, since its domain is a compact set. Then, it follows from Lebesgue's dominated convergence theorem that
\begin{align*}
&\int_0^1 \bigg( \int_{-T}^0 D_x^{\text{ac}}\Uc(\eta + \lambda(\eta_{1,\delta}-\eta))(\eta(x)-\eta_{1,\delta}(x)) dx\bigg) d\lambda \\
&= \ \int_0^1 \bigg( \int_{-T}^0 f(x,\lambda,\delta) dx\bigg) d\lambda \ \overset{\delta\rightarrow0^+}{\longrightarrow} \ \int_0^1 \bigg( \int_{-T}^0 f(x,\lambda,0) dx\bigg) d\lambda \\
&= \ \int_0^1 \bigg( \int_{-T}^0 D_x^{\text{ac}}u(\eta + \lambda(\eta_1-\eta))(\eta(x)-\eta_1(x)) dx\bigg) d\lambda,
\end{align*}
from which we deduce \eqref{E:u-u1}, since $\Uc(\eta_{1,\delta})\rightarrow u(\eta_1)$ as $\delta\rightarrow0^+$. Taking $\eta_1(\cdot)=\eta(\cdot-\eps)1_{[-T,0[}+\eta(0)1_{\{0\}}$, we obtain
\begin{align*}
&\frac{u(\eta) - u(\eta(\cdot-\eps)1_{[-T,0[}+\eta(0)1_{\{0\}})}{\eps} \\
&= \int_0^1 \bigg( \int_{-T}^0 D_x^{\text{ac}} u\big(\eta + \lambda \big(\eta(\cdot-\eps)-\eta(\cdot)\big) 1_{[-T,0[}\big) \frac{\eta(x)-\eta(x-\eps)}{\eps} dx \bigg) d\lambda \\
&= \ I_1(\eta,\eps) + I_2(\eta,\eps)+ I_3(\eta,\eps),
\end{align*}
where
\begin{align*}
I_1(\eta,\eps) \ &:= \ \int_0^1 \bigg(\int_{-T}^0 \eta(x)\frac{1}{\eps} \Big( D_x^{\text{ac}} u\big(\eta + \lambda \big(\eta(\cdot-\eps)-\eta(\cdot)\big) 1_{[-T,0[}\big) \\
&\quad \ - D_{x+\eps}^{\text{ac}} u\big(\eta + \lambda \big(\eta(\cdot-\eps)-\eta(\cdot)\big) 1_{[-T,0[}\big) \Big) dx \bigg) d\lambda, \\
I_2(\eta,\eps) \ &:= \ \frac{1}{\eps}\int_0^1 \bigg( \int_{-\eps}^0 \eta(x) D_{x+\eps}^{\text{ac}} u\big(\eta + \lambda \big(\eta(\cdot-\eps)-\eta(\cdot)\big) 1_{[-T,0[}\big) dx \bigg) d\lambda, \\
I_3(\eta,\eps) \ &:= \ - \frac{1}{\eps} \int_0^1 \bigg( \int_{-T-\eps}^{-T} \eta(x) D_{x+\eps}^{\text{ac}} u\big(\eta + \lambda \big(\eta(\cdot-\eps)-\eta(\cdot)\big) 1_{[-T,0[}\big) dx \bigg) d\lambda.
\end{align*}
Notice that, since $\eta(x)=0$ for $x<-T$, we see that $I_2(\eta,\eps)=0$. Moreover, since $D_x^{\text{ac}} u(\cdot)=D_0^{\text{ac}} u(\cdot)$, for $x\geq 0$, and $\eta + \lambda (\eta(\cdot-\eps)-\eta(\cdot)) 1_{[-T,0[}\rightarrow\eta$ in $\mathscr C([-T,0])$ as $\eps\rightarrow0^+$, it follows that (using the continuity of $D_x^{\text{ac}}u$ from $\mathscr C([-T,0])$ into $BV([-T,0])$, which implies that $D_0^{\text{ac}}u(\eta + \lambda (\eta(\cdot-\eps)-\eta(\cdot)) 1_{[-T,0[})\rightarrow D_0^{\text{ac}}u(\eta)$ as $\eps\rightarrow0^+$)
\begin{align*}
&\frac{1}{\eps}\int_{-\eps}^0 \eta(x)D_{x+\eps}^{\text{ac}} u\big(\eta + \lambda \big(\eta(\cdot-\eps)-\eta(\cdot)\big) 1_{[-T,0[}\big) dx \notag \\
&= \ \frac{1}{\eps}\int_{-\eps}^0 \eta(x) dx\,D_0^{\text{ac}}u\big(\eta + \lambda \big(\eta(\cdot-\eps)-\eta(\cdot)\big) 1_{[-T,0[}\big) \ \overset{\eps\rightarrow0^+}{\longrightarrow} \ \eta(0)D_0^{\text{ac}}u(\eta).
\end{align*}
Finally, concerning $I_1(\eta,\eps)$, from Fubini's theorem we obtain (denoting $\eta_{\eps,\lambda} := \eta + \lambda (\eta(\cdot-\eps)-\eta(\cdot)) 1_{[-T,0[}$)
\begin{align*}
I_1(\eta,\eps) \ &= \ \int_0^1 \bigg(\int_{-T}^0 \eta(x)\frac{1}{\eps} \Big( D_x^{\text{ac}} u(\eta_{\eps,\lambda}) - D_{x+\eps}^{\text{ac}} u(\eta_{\eps,\lambda}) \Big) dx \bigg) d\lambda \\
&= \ -\int_0^1 \bigg(\int_{-T}^0 \eta(x)\frac{1}{\eps} \bigg( \int_{]x,x+\eps]} D_{dy}^{\text{ac}} u(\eta_{\eps,\lambda})\bigg) dx \bigg) d\lambda \\
&= \ -\int_0^1 \bigg(\int_{]-T,\eps]} \frac{1}{\eps} \bigg( \int_{(-T)\vee(y-\eps)}^{0\wedge y} \eta(x) dx\bigg) D_{dy}^{\text{ac}} u(\eta_{\eps,\lambda}) \bigg) d\lambda \ = \ I_{11}(\eta,\eps) + I_{12}(\eta,\eps),
\end{align*}
where
\begin{align*}
I_{11}(\eta,\eps) \ &:= \ -\int_0^1 \bigg(\int_{]-T,\eps]} \frac{1}{\eps} \bigg( \int_{(-T)\vee(y-\eps)}^{0\wedge y} \eta(x) dx\bigg) \Big(D_{dy}^{\text{ac}} u(\eta_{\eps,\lambda}) - D_{dy}^{\text{ac}} u(\eta)\Big) \bigg) d\lambda, \\
I_{12}(\eta,\eps) \ &:= \ -\int_0^1 \bigg(\int_{]-T,\eps]} \frac{1}{\eps} \bigg( \int_{(-T)\vee(y-\eps)}^{0\wedge y} \eta(x) dx\bigg) D_{dy}^{\text{ac}} u(\eta) \bigg) d\lambda \\
&= \ -\bigg(\int_{]-T,\eps]} \frac{1}{\eps} \bigg( \int_{(-T)\vee(y-\eps)}^{0\wedge y} \eta(x) dx\bigg) D_{dy}^{\text{ac}} u(\eta).
\end{align*}
Recalling that $D_x^{\text{ac}} u(\cdot)=D_0^{\text{ac}} u(\cdot)$, for $x\geq 0$, we see that in $I_{11}(\eta,\eps)$ and $I_{12}(\eta,\eps)$ the integrals on $]-T,\eps]$ are equal to the same integrals on $]-T,0]$, i.e.,
\begin{align*}
I_{11}(\eta,\eps) \ &= \ -\int_0^1 \bigg(\int_{]-T,0]} \frac{1}{\eps} \bigg( \int_{(-T)\vee(y-\eps)}^{0\wedge y} \eta(x) dx\bigg) \Big(D_{dy}^{\text{ac}} u(\eta_{\eps,\lambda}) - D_{dy}^{\text{ac}} u(\eta)\Big) \bigg) d\lambda \\
&= \ -\int_0^1 \bigg(\int_{]-T,0]} \frac{1}{\eps} \bigg( \int_{y-\eps}^y \eta(x) dx\bigg) \Big(D_{dy}^{\text{ac}} u(\eta_{\eps,\lambda}) - D_{dy}^{\text{ac}} u(\eta)\Big) \bigg) d\lambda, \\
I_{12}(\eta,\eps) \ &= \ -\int_{]-T,0]} \frac{1}{\eps} \bigg( \int_{(-T)\vee(y-\eps)}^{0\wedge y} \eta(x) dx\bigg) D_{dy}^{\text{ac}} u(\eta) \ = \ -\int_{]-T,0]} \frac{1}{\eps} \bigg( \int_{y-\eps}^y \eta(x) dx\bigg) D_{dy}^{\text{ac}} u(\eta).
\end{align*}
Now, observe that
\[
|I_{11}(\eta,\eps)| \ \leq \ \|\eta\|_\infty \|D_\cdot^\text{ac}u(\eta_{\eps,\lambda})-D_\cdot^\text{ac}u(\eta)\|_{\text{Var}([-T,0])} \ \overset{\eps\rightarrow0^+}{\longrightarrow} \ 0.
\]
Moreover, since $\eta$ is continuous at $y\in]-T,0]$, we deduce that $\int_{y-\eps}^y \eta(x) dx/\eps\rightarrow\eta(y)$ as $\eps\rightarrow0^+$. Therefore, by Lebesgue's dominated convergence theorem, we get
\[
I_{12}(\eta,\eps) \ \overset{\eps\rightarrow0^+}{\longrightarrow} \ -\int_{]-T,0]} \eta(y) D_{dy}^{\text{ac}} u(\eta).
\]
In conclusion, we have
\[
D^H\Uc(\eta) \ = \ \eta(0)D_0^{\text{ac}}u(\eta) - \int_{]-T,0]} \eta(y) D_{dy}^{\text{ac}} u(\eta),
\]
which gives \eqref{E:DH=Dacdeta} using the integration by parts formula \eqref{E:IbyP} and noting that we can suppose, without loss of generality, $D_{0^-}^{\text{ac}}\Uc(\eta)=D_0^{\text{ac}}\Uc(\eta)$.
\ep

\vspace{3mm}

For our second representation result of $D^H\Uc$ we need the following generalization of the deterministic backward integral when the integrand is a measure.

\begin{Definition}
\label{D:DeterministicIntegralMeasure}
Let $f\colon[-T,0]\rightarrow\R$ be a c\`{a}dl\`{a}g function and $g\in\Mc([-T,0])$. Suppose that the following limit
\[
\int_{[-T,0]}g(ds)d^+f(s) \ := \ \lim_{\eps\rightarrow0^+}\int_{[-T,0]} g(ds)\frac{f_{\overline J}(s)-f_{\overline J}(s-\eps)}{\eps},
\]
exists and it is finite. Then, the obtained quantity is denoted by $\int_{[-T,0]} gd^+f$ and called \textbf{$($determin\-istic, definite$)$ backward integral of $g$ with respect to $f$ $($on $[-T,0]$}$)$.
\end{Definition}

\begin{Remark}
{\rm
Notice that if $g$ is absolutely continuous with density c\`adl\`ag $($still denoted by $g$$)$ then Definition 
\ref{D:DeterministicIntegralMeasure} is compatible
with the one in Definition  \ref{D:DeterministicIntegral}.
\ep
}
\end{Remark}

\begin{Proposition}
\label{P:DH_SecondOrder}
Let $\Uc\colon C([-T,0])\rightarrow\R$ be twice continuously Fr\'echet differentiable such that
\[
D^2\Uc\colon C([-T,0]) \ \longrightarrow \ \chi_0\subset(C([-T,0])\hat\otimes_\pi C([-T,0]))^*\text{ continuously with respect to $\chi_0$.}
\]
Assume that there exist continuous extensions (necessarily unique)
\[
u\colon\mathscr C([-T,0])\rightarrow\R, \qquad\qquad D_{dx\,dy}^2u\colon\mathscr C([-T,0])\rightarrow\chi_0
\]
of $\Uc$ and $D_{dx\,dy}^2\Uc$, respectively.
Let $\eta\in C([-T,0])$ be such that the (deterministic) 
quadratic variation on $[-T,0]$ exists and suppose also the following. 
\begin{enumerate}
\item[\textup{(i)}] $D_x^{2,Diag}\Uc(\eta)$, the diagonal element of the second-order derivative at $\eta$, has a set of discontinuity which has null measure with respect to $[\eta]$ $($in particular, if it is countable$)$.
\item[\textup{(ii)}] The horizontal derivative $D^H\Uc(\eta)$ exists at $\eta$.
\end{enumerate}
Then
\begin{equation}
\label{E:DH=SecondOrder}
D^H \Uc(\eta) \ = \ \int_{[-T,0]} D_{dx}^\perp \Uc(\eta) d^+ \eta(x) - \frac{1}{2}\int_{[-T,0]} D_x^{2,Diag}\Uc(\eta) d[\eta](x).
\end{equation}
In particular, the backward integral in \eqref{E:DH=SecondOrder} exists.
\end{Proposition}
\textbf{Proof.} Let $\eta$ be as in the statement of the proposition. Then, using the definition of $D^H\Uc(\eta)$ we are led to consider the following increment for the function $u$:
\begin{equation}
\label{E:ProofDH}
\frac{u(\eta) - u(\eta(\cdot-\eps)1_{[-T,0[}+\eta(0)1_{\{0\}})}{\eps},
\end{equation}
with $\eps>0$. Our aim is to expand \eqref{E:ProofDH} using some Taylor formula. To this end, we begin noting that, since $\Uc$ is $C^2$ Fr\'echet, for every $\eta_1\in C([-T,0])$ the following standard Taylor's expansion holds:
\begin{align*}
&\Uc(\eta_1) \ = \ \Uc(\eta) + \int_{[-T,0]}\!\! D_{dx}\Uc(\eta) \big(\eta_1(x) - \eta(x)\big) + \frac{1}{2}\int_{[-T,0]^2} \!\! D_{dx\,dy}^2\Uc(\eta) \big(\eta_1(x)-\eta(x)\big) \big(\eta_1(y)-\eta(y)\big) \notag \\
& + \int_0^1 (1-\lambda)\bigg( \int_{[-T,0]^2} \Big( D_{dx\,dy}^2\Uc(\eta + \lambda (\eta_1-\eta)) - D_{dx\,dy}^2\Uc(\eta) \Big) \big(\eta_1(x)-\eta(x)\big) \big(\eta_1(y)-\eta(y)\big) \bigg) d\lambda. \notag
\end{align*}
Now, using the density of $C_{\eta(0)}([-T,0])$ into $\mathscr C_{\eta(0)}([-T,0])$ with respect to the topology of the space $\mathscr C([-T,0])$ and proceeding as in the proof of Proposition \ref{P:DH=Dacdeta}, we deduce the following Taylor's formula for $u$:
\begin{align}
\label{E:SecondOrder1}
&\frac{u(\eta) - u(\eta(\cdot-\eps)1_{[-T,0[}+\eta(0)1_{\{0\}})}{\eps} \ = \ \int_{[-T,0]} D_{dx}^\perp \Uc(\eta) \frac{\eta(x) - \eta(x-\eps)}{\eps} \\
& - \frac{1}{2}\int_{[-T,0]^2} D_{dx\,dy}^2 \Uc(\eta) \frac{(\eta(x)-\eta(x-\eps)) (\eta(y)-\eta(y-\eps))}{\eps} 1_{[-T,0[\times[-T,0[}(x,y) \notag \\
&\quad \ - \int_0^1 (1-\lambda)\bigg( \int_{[-T,0]^2} \Big( D_{dx\,dy}^2 u(\eta + \lambda (\eta(\cdot-\eps)-\eta(\cdot))1_{[-T,0[}) \notag \\
& - D_{dx\,dy}^2 \Uc(\eta) \Big) \frac{(\eta(x)-\eta(x-\eps)) (\eta(y)-\eta(y-\eps))}{\eps} 1_{[-T,0[\times[-T,0[}(x,y) \bigg) d\lambda. \notag
\end{align}
Recalling the definition of $\chi_0$ given in Remark \ref{R:Chi-subspace}, denoting by $D_{x\,y}^{2,L^2} \Uc(\eta)\in L^2([-T,0]^2)$ the element $g_1$ and by $D_{x}^{2,Diag}\Uc(\eta)\in L^\infty([-T,0])$ the diagonal element $g_4$ of $D_{dx\,dy}^2\Uc(\eta)$, we notice that (due to the presence of the indicator function $1_{[-T,0[\times[-T,0[}$)
\begin{align*}
&\int_{[-T,0]^2} D_{dx\,dy}^2 \Uc(\eta) \frac{(\eta(x)-\eta(x-\eps)) (\eta(y)-\eta(y-\eps))}{\eps} 1_{[-T,0[\times[-T,0[}(x,y) \\
&= \ \int_{[-T,0]^2} D_{x\,y}^{2,L^2} \Uc(\eta) \frac{(\eta(x)-\eta(x-\eps)) (\eta(y)-\eta(y-\eps))}{\eps} dx\,dy \\
&\quad \ + \int_{[-T,0]} D_x^{2,Diag} \Uc(\eta) \frac{(\eta(x)-\eta(x-\eps))^2}{\eps} dx.
\end{align*}
We denote by $D_{x\,y}^{2,L^2} u$ and $D_{x}^{2,Diag}u$ the extensions of $D_{x\,y}^{2,L^2}\Uc$ and $D_{x}^{2,Diag}\Uc$ to $\mathscr C([-T,0])$, respectively, which are continuous. In particular, \eqref{E:SecondOrder1} becomes
\begin{equation}
\label{E:SecondOrder2}
\frac{u(\eta) - u(\eta(\cdot-\eps)1_{[-T,0[}+\eta(0)1_{\{0\}})}{\eps} \ = \ I_1(\eps) + I_2(\eps) + I_3(\eps) + I_4(\eps) + I_5(\eps),
\end{equation}
where
\begin{align*}
I_1(\eps) \ &:= \ \int_{[-T,0]} D_{dx}^\perp \Uc(\eta) \frac{\eta(x) - \eta(x-\eps)}{\eps}, \\
I_2(\eps) \ &:= \ - \frac{1}{2}\int_{[-T,0]^2} D_{x\,y}^{2,L^2} \Uc(\eta) \frac{(\eta(x)-\eta(x-\eps)) (\eta(y)-\eta(y-\eps))}{\eps} dx\,dy, \\
I_3(\eps) \ &:= \ - \frac{1}{2}\int_{[-T,0]} D_x^{2,Diag} \Uc(\eta) \frac{(\eta(x)-\eta(x-\eps))^2}{\eps} dx, \\
I_4(\eps) \ &:= \ - \int_0^1 (1-\lambda)\bigg( \int_{[-T,0]^2} \Big( D_{x\,y}^{2,L^2} u(\eta + \lambda (\eta(\cdot-\eps)-\eta(\cdot))1_{[-T,0[}) \\
&\quad \ - D_{x\,y}^{2,L^2} \Uc(\eta) \Big) \frac{(\eta(x)-\eta(x-\eps)) (\eta(y)-\eta(y-\eps))}{\eps} dx\,dy \bigg) d\lambda, \\
I_5(\eps) \ &:= \ - \int_0^1 (1-\lambda)\bigg( \int_{[-T,0]} \Big( D_x^{2,Diag} u(\eta + \lambda (\eta(\cdot-\eps)-\eta(\cdot))1_{[-T,0[}) \\
&\quad \ - D_x^{2,Diag} \Uc(\eta) \Big) \frac{(\eta(x)-\eta(x-\eps))^2}{\eps} dx \bigg) d\lambda.
\end{align*}
Firstly, we shall prove that
\begin{equation}
\label{E:I2-->0}
I_2(\eps) \ \overset{\eps\rightarrow0^+}{\longrightarrow} \ 0.
\end{equation}
To this end, for every $\eps>0$, define the operator $T_\eps\colon L^2([-T,0]^2)\rightarrow\R$ as follows:
\[
T_\eps \, g \ = \ \int_{[-T,0]^2} g(x,y) \frac{(\eta(x)-\eta(x-\eps)) (\eta(y)-\eta(y-\eps))}{\eps} dx\,dy, \qquad \forall\, g\in L^2([-T,0]^2).
\]
Then $T_\eps\in L^2([-T,0])^*$. Indeed, from Cauchy-Schwarz inequality,
\begin{align*}
|T_\eps\,g| \ &\leq \ \|g\|_{L^2([-T,0]^2)} \sqrt{\int_{[-T,0]^2}\frac{(\eta(x)-\eta(x-\eps))^2 (\eta(y)-\eta(y-\eps))^2}{\eps^2}dx\,dy} \\
&= \ \|g\|_{L^2([-T,0]^2)} \int_{[-T,0]} \frac{(\eta(x)-\eta(x-\eps))^2}{\eps} dx
\end{align*}
and this last quantity is  bounded with respect to $\eps$ since the quadratic variation of $\eta$ on $[-T,0]$ exists. In particular, we have proved that for every $g\in L^2([-T,0]^2)$ there exists a constant $M_g\geq0$ such that
\[
\sup_{0 < \eps < 1} |T_\eps\,g| \ \leq \ M_g.
\]
It follows from Banach-Steinhaus theorem that there exists a constant $M\geq0$ such that
\begin{equation}
\label{E:BanachSteinhaus}
\sup_{0 < \eps < 1} \|T_\eps\|_{L^2([-T,0])^*} \ \leq \ M.
\end{equation}
Now, let us consider the set $\Sc := \{g\in L^2([-T,0]^2)\colon g(x,y) = e(x)f(y),\text{ with }e,f\in C^1([-T,0])\}$, which is dense in $L^2([-T,0]^2)$. Let us show that
\begin{equation}
\label{E:T_eps-->0}
T_\eps\,g \ \overset{\eps\rightarrow0^+}{\longrightarrow} \ 0, \qquad \forall\,g\in\Sc.
\end{equation}
Fix $g\in\Sc$, with $g(x,y)=e(x)f(y)$ for any $(x,y)\in[-T,0]$, then
\begin{equation}
\label{E:T_eps_Sc}
T_\eps\,g \ = \ \frac{1}{\eps} \int_{[-T,0]} e(x) \big(\eta(x) - \eta(x-\eps)\big) dx \int_{[-T,0]} f(y) \big(\eta(y) - \eta(y-\eps)\big) dy.
\end{equation}
We have
\begin{align*}
&\bigg|\int_{[-T,0]} e(x) \big(\eta(x) - \eta(x-\eps)\big) dx\bigg| \ = \ \bigg|\int_{[-T,0]} \big(e(x) - e(x+\eps)\big) \eta(x) dx \\
&- \int_{[-T-\eps,-T]}e(x+\eps)\eta(x) dx + \int_{[-\eps,0]} e(x+\eps)\eta(x) dx\bigg| \ \leq \ \eps \bigg(\int_{[-T,0]}|\dot e(x)|dx + 2\|e\|_\infty\bigg)\|\eta\|_\infty.
\end{align*}
Similarly,
\[
\bigg|\int_{[-T,0]} f(y) \big(\eta(y) - \eta(y-\eps)\big) dy\bigg| \ \leq \ \eps \bigg(\int_{[-T,0]}|\dot f(y)|dy + 2\|f\|_\infty\bigg)\|\eta\|_\infty.
\]
Therefore, from \eqref{E:T_eps_Sc} we find
\[
|T_\eps\,g| \ \leq \ \eps \bigg(\int_{[-T,0]}|\dot e(x)|dx + 2\|e\|_\infty\bigg) \bigg(\int_{[-T,0]}|\dot f(y)|dy + 2\|f\|_\infty\bigg) \|\eta\|_\infty^2,
\]
which converges to zero as $\eps$ goes to zero and therefore \eqref{E:T_eps-->0} is established. This in turn implies that
\begin{equation}
\label{E:T_eps-->0bis}
T_\eps\,g \ \overset{\eps\rightarrow0^+}{\longrightarrow} \ 0, \qquad \forall\,g\in L^2([-T,0]^2).
\end{equation}
Indeed, fix $g\in L^2([-T,0]^2)$ and let $(g_n)_n\subset\Sc$ be such that $g_n\rightarrow g$ in $L^2([-T,0]^2)$. Then
\[
|T_\eps\,g| \ \leq \ |T_\eps(g-g_n)| + |T_\eps\,g_n| \ \leq \ \|T_\eps\|_{L^2([-T,0]^2)^*} \|g-g_n\|_{L^2([-T,0]^2)} + |T_\eps\,g_n|.
\]
From \eqref{E:BanachSteinhaus} it follows that
\[
|T_\eps\,g| \ \leq \ M \|g-g_n\|_{L^2([-T,0]^2)} + |T_\eps\,g_n|,
\]
which implies $\limsup_{\eps\rightarrow0^+}|T_\eps\,g| \leq M \|g-g_n\|_{L^2([-T,0]^2)}$. Sending $n$ to infinity, we deduce \eqref{E:T_eps-->0bis} and finally \eqref{E:I2-->0}.

Let us now consider the term $I_3(\eps)$ in \eqref{E:SecondOrder2}. Since the quadratic variation $[\eta]$ exists, it follows from Portmanteau's theorem and item (i) that
\[
I_3(\eps) \ = \ \int_{[-T,0]} D_x^{2,Diag} \Uc(\eta) \frac{(\eta(x)-\eta(x-\eps))^2}{\eps} dx \ \underset{\eps\rightarrow0^+}{\longrightarrow} \ \int_{[-T,0]} D_x^{2,Diag}\Uc(\eta) d[\eta](x).
\]
Regarding the term $I_4(\eps)$ in \eqref{E:SecondOrder2}, let $\phi_\eta\colon[0,1]^2\rightarrow L^2([-T,0]^2)$ be given by
\[
\phi_\eta(\eps,\lambda)(\cdot,\cdot) \ = \ D_{\cdot\,\cdot}^{2,L^2} u\big(\eta + \lambda (\eta(\cdot-\eps)-\eta(\cdot))1_{[-T,0[}\big).
\]
By assumption, $\phi_\eta$ is a continuous map, and hence it is uniformly continuous, since $[0,1]^2$ is a compact set. Let $\rho_{\phi_\eta}$ denote the continuity modulus of $\phi_\eta$, then
\begin{align*}
&\big\|D_{\cdot\,\cdot}^{2,L^2} u\big(\eta + \lambda (\eta(\cdot-\eps)-\eta(\cdot))1_{[-T,0[}\big) - D_{\cdot\,\cdot}^{2,L^2} \Uc(\eta)\big\|_{L^2([-T,0]^2)} \\
&= \ \|\phi_\eta(\eps,\lambda) - \phi_\eta(0,\lambda)\|_{L^2([-T,0]^2)} \ \leq \ \rho_{\phi_\eta}(\eps).
\end{align*}
This implies, by Cauchy-Schwarz inequality,
\begin{align*}
&\bigg|\int_0^1 (1-\lambda)\bigg( \int_{[-T,0]^2} \Big( D_{x\,y}^{2,L^2} u(\eta + \lambda (\eta(\cdot-\eps)-\eta(\cdot))1_{[-T,0[}) \\
&- D_{x\,y}^{2,L^2} \Uc(\eta) \Big) \frac{(\eta(x)-\eta(x-\eps)) (\eta(y)-\eta(y-\eps))}{\eps} dx\,dy \bigg) d\lambda\bigg| \\
&\leq \ \int_0^1(1-\lambda)\big\|D_{\cdot\,\cdot}^{2,L^2}u(\eta+\lambda(\eta(\cdot-\eps)-\eta(\cdot))1_{[-T,0]}) \\
&- D_{\cdot\,\cdot}^{2,L^2}\Uc(\eta)\big\|_{L^2([-T,0]^2)}\sqrt{\int_{[-T,0]^2}\frac{(\eta(x)-\eta(x-\eps))^2 (\eta(y)-\eta(y-\eps))^2}{\eps^2}dx\,dy}\,d\lambda \\
&\leq \ \int_0^1(1-\lambda)\rho_{\phi_\eta}(\eps) \bigg(\int_{[-T,0]} \frac{(\eta(x)-\eta(x-\eps))^2}{\eps} dx \bigg)d\lambda \\
&= \ \frac{1}{2}\rho_{\phi_\eta}(\eps) \int_{[-T,0]} \frac{(\eta(x)-\eta(x-\eps))^2}{\eps} dx \ \overset{\eps\rightarrow0^+}{\longrightarrow} \ 0.
\end{align*}
Finally, we consider the term $I_5(\eps)$ in \eqref{E:SecondOrder2}. Define $\psi_\eta\colon[0,1]^2\rightarrow L^\infty([-T,0])$ as follows:
\[
\psi_\eta(\eps,\lambda)(\cdot) \ = \ D_\cdot^{2,Diag} u\big(\eta + \lambda (\eta(\cdot-\eps)-\eta(\cdot))1_{[-T,0[}\big).
\]
We see that $\psi_\eta$ is uniformly continuous. Let $\rho_{\psi_\eta}$ denote the continuity modulus of $\psi_\eta$, then
\begin{align*}
&\big\|D_\cdot^{2,Diag} u\big(\eta + \lambda (\eta(\cdot-\eps)-\eta(\cdot))1_{[-T,0[}\big) - D_\cdot^{2,Diag} \Uc(\eta)\big\|_{L^\infty([-T,0])} \\
&= \ \|\psi_\eta(\eps,\lambda) - \psi_\eta(0,\lambda)\|_{L^\infty([-T,0])} \ \leq \ \rho_{\psi_\eta}(\eps).
\end{align*}
Therefore, we have
\begin{align*}
&\bigg|\int_0^1 (1-\lambda)\bigg( \int_{[-T,0]} \Big( D_x^{2,Diag} u(\eta + \lambda (\eta(\cdot-\eps)-\eta(\cdot))1_{[-T,0[}) \\
&- D_x^{2,Diag} \Uc(\eta) \Big) \frac{(\eta(x)-\eta(x-\eps))^2}{\eps} dx \bigg) d\lambda\bigg| \\
&\leq \ \int_0^1(1-\lambda)\bigg(\int_{[-T,0]}\rho_{\psi_\eta}(\eps) \frac{(\eta(x)-\eta(x-\eps))^2}{\eps} dx \bigg) d\lambda \\
&= \ \frac{1}{2}\rho_{\psi_\eta}(\eps)\int_{[-T,0]} \frac{(\eta(x)-\eta(x-\eps))^2}{\eps} dx \ \overset{\eps\rightarrow0^+}{\longrightarrow} \ 0.
\end{align*}
In conclusion, we have proved that all the integral terms in the right-hand side of \eqref{E:SecondOrder2}, unless $I_1(\eps)$, admit a limit when $\eps$ goes to zero. Since the left-hand side admits a limit, namely $D^H\Uc(\eta)$, we deduce that the backward integral
\[
I_1(\eps) \ = \ \int_{[-T,0]} D_{dx}^\perp \Uc(\eta) \frac{\eta(x) - \eta(x-\eps)}{\eps} \ \overset{\eps\rightarrow0^+}{\longrightarrow} \ \int_{[-T,0]} D_{dx}^\perp \Uc(\eta) d^+\eta(x)
\]
exists and it is finite, which concludes the proof.
\ep

\section{Path-dependent SDE and Kolmogorov equation}

\label{S3}

\subsection{The framework}

Consider on $(\Omega,\Fc,\P)$ a real Brownian motion $W=(W_t)_{t\geq0}$. We denote by $\F=(\Fc_t)_{t\geq0}$ the natural filtration generated by $W$, completed with the $\P$-null sets of $\Fc$. Fix a finite time horizon $T\in]0,\infty[$ and let $C([-T,0])$ be the Banach space of all continuous paths $\eta\colon[-T,0]\rightarrow\R$ endowed with the supremum norm $\|\eta\|_\infty=\sup_{x\in[-T,0]}|\eta(x)|$. For any $(t,\eta)\in[0,T]\times C([-T,0])$ consider the following \emph{path-dependent SDE} (recall that $\X$ is the window process associated to $X$, see \eqref{WindowProcess})
\begin{equation}
\label{SDE}
\begin{cases}
dX_s = \ b(s,\mathbb X_s)dt + \sigma(s,\mathbb X_s)dW_s, \qquad\qquad & t\leq s\leq T, \\
X_s \ = \ \eta(s-t), & -T+t\leq s\leq t,
\end{cases}
\end{equation}
where the coefficients $b$ and $\sigma$ satisfy the following assumption.

\begin{itemize}
\item[\textbf{(H1)}] $b,\sigma\colon[0,T]\times C([-T,0])\rightarrow\R$ are Borel measurable functions satisfying, for all $t\in[0,T]$, $\eta,\eta'\in C([-T,0])$,
\begin{align*}
|b(t,\eta)| + |\sigma(t,\eta)| \ &\leq \ M_1(1 + \|\eta\|_\infty), \\
|b(t,\eta) - b(t,\eta')| + |\sigma(t,\eta) - \sigma(t,\eta')| \ &\leq \ L_1\|\eta - \eta'\|_\infty,
\end{align*}
for some positive constants $M_1$ and $L_1$.
\end{itemize}

\begin{Proposition}
Under Assumption {\bf (H1)}, for any $(t,\eta)\in[0,T]\times C([-T,0])$ there exists a unique (up to indistinguishability) $\F$-adapted continuous process $X^{t,\eta}=(X_s^{t,\eta})_{[-T+t,T]}$ strong solution to equation \eqref{SDE}. Moreover, for any $p\geq2$ there exists a positive constant $C_p$ (depending only on $p,T,M_1$) such that
\begin{equation}
\label{EstimateX}
\E\Big[\sup_{s\in[-T+t,T]}\big|X_s^{t,\eta}\big|^p\Big] \ \leq \ C_p \big(1 + \|\eta\|_\infty^p\big).
\end{equation}
\end{Proposition}
\textbf{Proof.}
Since Hypotheses (14.15) and (14.22) in \cite{jacod79} are satisfied under \textbf{(H1)}, the existence and uniqueness part follows from Theorem 14.23 in \cite{jacod79}. Concerning \eqref{EstimateX}, raising to the $p$-th power both sides of equation \eqref{SDE}, recalling that $(x_1+\cdots+x_n)^p\leq n^{p-1}(x_1^p+\cdots+x_n^p),$ for any $x_1,\ldots,x_n\geq0$, we obtain ($\X^{t,\eta}$ is the window process associated to $X^{t,\eta}$, see \eqref{WindowProcess})
\[
|X_s^{t,\eta}|^p \ \leq \ 3^{p-1}\bigg\{|\eta(0)|^p + \bigg|\int_t^s b(r,\mathbb X_r^{t,\eta}) dr\bigg|^p + \bigg|\int_t^s \sigma(r,\mathbb X_r^{t,\eta}) dW_r\bigg|^p\bigg\}.
\]
Set $\Xc_s^{t,\eta}=\sup_{u\in[-T+t,s]}|X_u^{t,\eta}|$, for all $s\in[t,T]$. Then, we have
\begin{equation}
\label{EstimateX_Proof1}
|\Xc_s^{t,\eta}|^p \ \leq \ 3^{p-1}\bigg\{\|\eta\|_\infty^p + \sup_{u\in[t,s]}\bigg|\int_t^u b(r,\mathbb X_r^{t,\eta}) dr\bigg|^p + \sup_{u\in[t,s]}\bigg|\int_t^u \sigma(r,\mathbb X_r^{t,\eta}) dW_r\bigg|^p\bigg\}.
\end{equation}
From H\"older's inequality, we get
\begin{equation}
\label{EstimateX_Proof2}
\sup_{u\in[t,s]}\bigg|\int_t^u b(r,\mathbb X_r^{t,\eta}) dr\bigg|^p \ \leq \ (s-t)^{p-1} \int_t^s \big| b(r,\mathbb X_r^{t,\eta}) \big|^p dr.
\end{equation}
On the other hand, from Burkholder-Davis-Gundy's inequality and H\"older's inequality, it follows that there exists a positive constant $c_p$, depending only on $p$, such that
\begin{align}
\label{EstimateX_Proof3}
\E\bigg[\sup_{u\in[t,s]}\bigg|\int_t^u \sigma(r,\mathbb X_r^{t,\eta}) dW_r\bigg|^p\bigg] \ &\leq \ c_p \E\bigg[\bigg|\int_t^s \big|\sigma(r,\mathbb X_r^{t,\eta})\big|^2 dr\bigg|^\frac{p}{2}\bigg] \notag \\
&\leq \ c_p (s-t)^{\frac{p}{2}-1} \E\bigg[\int_t^s \big|\sigma(r,\mathbb X_r^{t,\eta})\big|^p dr\bigg].
\end{align}
Taking the expectation in \eqref{EstimateX_Proof1}, then exploiting \eqref{EstimateX_Proof2} and \eqref{EstimateX_Proof3}, afterwards using the linear growth condition of $b$ and $\sigma$ in {\bf (H1)}, we see that there exists a positive constant $\bar C_p$, depending only on $p,T,M_1$, such that
\[
|\Xc_s^{t,\eta}|^p \ \leq \ \bar C_p \Big(1 + \|\eta\|_\infty^p + \int_t^s \E[|\Xc_r^{t,\eta}|^p] dr\Big), \qquad \forall\,s\in[t,T].
\]
We can now deduce estimate \eqref{EstimateX} from an application of Gronwall's inequality.
\ep

\vspace{3mm}

Our aim is to study the following \emph{path-dependent Kolmogorov equation}, which turns out to be related to \eqref{SDE}:
\begin{equation}
\label{KolmEq}
\begin{cases}
\partial_t\Uc + D^H\Uc + b(t,\eta)D^V\Uc \\
+ \frac{1}{2}\sigma(t,\eta)^2D^{VV}\Uc + F(t,\eta,\Uc,\sigma(t,\eta)D^V\Uc) \ = \ 0, \;\;\; &\forall\,(t,\eta)\in[0,T[\times C([-T,0]), \\
\Uc(T,\eta) \ = \ H(\eta), &\forall\,\eta\in C([-T,0]),
\end{cases}
\end{equation}
where $F\colon[0,T]\times C([-T,0])\times\R\times\R\rightarrow\R$ and $H\colon C([-T,0])\rightarrow\R$ are Borel measurable functions.

\begin{Remark} \label{RNonAntic}
{\rm
As already recalled, our functionals are defined, differently from \cite{contfournie10}, in such a way that time $t$ and path $\eta$ are not related to each other, so that the non-anticipative property imposed by \cite{contfournie10} is automatically satisfied. More precisely, the non-anticipative property is implicit in the definition of our functionals, since given a pair $(t,\eta)$ the path $\eta$ always represents the past up to time $t$. Indeed, in general $\eta$ stands for a path of $\X_t$, which is the path of the process $X$ on $[-T+t,t]$.
\ep
}
\end{Remark}

\subsection{Strict solutions}

We introduce the concept of strict solution to equation \eqref{KolmEq} and then study its well-posedness.

\begin{Definition}
A map $\Uc\colon[0,T]\times C([-T,0])\rightarrow\R$ in $C^{1,2}(([0,T[\times\textup{past})\times\textup{present})\cap C([0,T]\times C([-T,0]))$, satisfying equation \eqref{KolmEq}, is called a \textbf{strict solution} to equation \eqref{KolmEq}.
\end{Definition}

Let us begin focusing on the uniqueness of equation \eqref{KolmEq}. Actually, we shall prove a stronger result, namely, that any strict solution of \eqref{KolmEq} can be represented in terms of a backward stochastic differential equation. For this probabilistic representation formula, it is convenient to introduce the following spaces of stochastic processes.

\begin{itemize}
\item $\S^p(t,T)$, $p\geq1$, $0 \leq t \leq T$, the set  of real-valued continuous $\F$-adapted stochastic processes $Y=(Y_s)_{t\leq s\leq T}$ satisfying
\[
\|Y\|_{_{\S^p(t,T)}}^p := \ \E\Big[ \sup_{t\leq s\leq T} |Y_s|^p \Big] \ < \ \infty.
\]
\item $\H^p(t,T)$, $p$ $\geq$ $1$, $0 \leq t \leq T$, the set of real-valued $\F$-predictable stochastic processes $Z=(Z_s)_{t\leq s\leq T}$ satisfying
\[
\|Z\|_{_{\H^p(t,T)}}^p := \ \E\bigg[\bigg(\int_t^T |Z_s|^2 ds\bigg)^{\frac{p}{2}}\bigg] \ < \ \infty.
\]
\end{itemize}

\begin{Theorem}
\label{T:UniqStrict}
Suppose that Assumption {\bf (H1)} holds. Let $F\colon[0,T]\times C([-T,0])\times\R\times\R\rightarrow\R$, $H\colon C([-T,0])\rightarrow\R$ be Borel measurable maps and consider a strict solution $\Uc\colon[0,T]\times C([-T,0])\rightarrow\R$ to equation \eqref{KolmEq}, satisfying, for some constants $C,m\geq0$,
\begin{align}
\label{PolGrowth}
|F(t,\eta,y,z) - F(t,\eta,y',z')| \ &\leq \ C\big(|y-y'| + |z-z'|\big), \notag \\
|H(\eta)| + |F(t,\eta,y,z)| + |\Uc(t,\eta)| \ &\leq \ C\big(1 + \|\eta\|_\infty^m\big),
\end{align}
for all $(t,\eta)\in[0,T]\times C([-T,0])$, $y,y',z,z'\in\R$. Then, we have
\begin{equation}
\label{U=Y}
\Uc(t,\eta) \ = \ Y_t^{t,\eta}, \qquad \forall\,(t,\eta)\in[0,T]\times C([-T,0]),
\end{equation}
where $(Y_s^{t,\eta},Z_s^{t,\eta})_{s\in[t,T]} = (\Uc(s,\mathbb X_s^{t,\eta}),(\sigma D^V\Uc)(s,\mathbb X_s^{t,\eta})1_{[t,T[}(s))_{s\in[t,T]}\in\S^2(t,T)\times\H^2(t,T)$ solves the backward SDE: $\P$-a.s.,
\begin{equation}
\label{BSDE}
Y_s^{t,\eta} \ = \ H(\mathbb X_T^{t,\eta}) + \int_s^T F(r,\mathbb X_r^{t,\eta},Y_r^{t,\eta},Z_r^{t,\eta}) dr - \int_s^T Z_r^{t,\eta} dW_r, \qquad t \leq s \leq T.
\end{equation}
As a consequence, there exists at most one strict solution to equation \eqref{KolmEq} satisfying a polynomial growth condition as in \eqref{PolGrowth}.
\end{Theorem}
\textbf{Proof.}
Take $(t,\eta)\in[0,T[\times C([-T,0])$ and define, for any $t \leq s \leq T$,
\[
Y_s^{t,\eta} \ = \ \Uc(s,\mathbb X_s^{t,\eta}), \qquad Z_s^{t,\eta} \ = \ \sigma(s,\mathbb X_s^{t,\eta}) D^V\Uc(s,\mathbb X_s^{t,\eta})1_{[t,T[}(s).
\]
Let $T_0\in[t,T[$, then an application of the functional It\^o formula \eqref{E:ItoTime} to $\Uc(s,\mathbb X_s^{t,\eta})$, yields (recall that $\Uc$ solves equation \eqref{KolmEq})
\begin{equation}
\label{E:BSDEIto}
Y_s^{t,\eta} \ = \ Y_{T_0}^{t,\eta} + \int_s^{T_0} F(r,\mathbb X_r^{t,\eta},Y_r^{t,\eta},Z_r^{t,\eta}) dr - \int_s^{T_0} Z_r^{t,\eta} dW_r, \qquad t \leq s \leq T_0.
\end{equation}
To conclude to the validity of \eqref{U=Y}  we have to take the limit $T_0\rightarrow T$ in 
\eqref{E:BSDEIto}. To this end, it is enough to have a uniform (with respect to $T_0$) bound on the norm of $Z^{t,\eta}$ in $\H^2(t,T_0)$. This follows from the following standard estimate for backward SDE, see  Proposition B.1 of Appendix in \cite{cosso_russo15b}, with $K = 0$:
\[
\E\int_t^{T_0} |Z_s^{t,\eta}|^2 ds \ \leq \ C \bigg(\|Y^{t,\eta}\|_{\S^2(t,T)}^2 + \E\int_t^T |F(s,\X_s^{t,\eta},0,0)|^2 ds\bigg), \qquad \forall\,T_0\in[t,T[,
\]
where $\bar C$ is a positive constant, depending only on $T$ and $C$. From estimate \eqref{EstimateX}, we get, for any $p\geq2$,
\begin{equation}
\label{EstimateSupX}
\E\Big[\sup_{t\leq s\leq T}\|\mathbb X_s^{t,\eta}\|_\infty^p\Big] \ < \ \infty.
\end{equation}
Notice that from the polynomial growth condition of $\Uc$ and \eqref{EstimateSupX} we have $\|Y^{t,\eta}\|_{\S^2(t,T)}<\infty$, therefore $Y\in\S^2(t,T)$. As a consequence, using monotone convergence theorem, we obtain
\[
\E\int_t^T |Z_s^{t,\eta}|^2 ds \ \leq \ C \bigg(\|Y^{t,\eta}\|_{\S^2(t,T)}^2 + \E\int_t^T |F(s,\mathbb X_s^{t,\eta},0,0)|^2 ds\bigg).
\]
This implies, using estimate \eqref{EstimateSupX} and the polynomial growth condition of $F$, that $Z\in\H^2(t,T)$. Moreover, from the uniform Lipschitz property of $F$ with respect to $(y,z)$, we see that $\E\int_t^T |F(r,$ $\mathbb X_r^{t,\eta},Y_r^{t,\eta},Z_r^{t,\eta})|^2 dr<\infty$. In conclusion, we can send $T_0\rightarrow T$ in \eqref{E:BSDEIto} and obtain the probabilistic representation formula \eqref{U=Y}.

As it is well-known that there exists a unique solution $(Y^{t,\eta},Z^{t,\eta})\in\S^2(t,T)\times\H^2(t,T)$ to the backward SDE \eqref{BSDE} (see, e.g., Theorem 4.1 in \cite{parpen90}), we deduce the uniqueness result for $\Uc$.
\ep

\vspace{3mm}

We now provide two existence results, i.e., Theorem \ref{T:ExistenceStrict} and Theorem \ref{T:ExistenceStrictPardouxPeng}, for strict solutions to equation \eqref{KolmEq}, when the coefficients have a \emph{cylindrical} form. First, we consider the case where $F$ does not depend on $(y,z)$, then we address the general semilinear case.

\begin{Theorem}
\label{T:ExistenceStrict}
Suppose that $F=F(t,\eta)$ does not depend on $(y,z)$ and there exists $N\in\N\backslash\{0\}$ such that, for all $(t,\eta)\in[0,T]\times C([-T,0])$,
\begin{align*}
b(t,\eta) \ &= \ \bar b\bigg(t,\int_{[-t,0]}\varphi_1(x+t)d^-\eta(x),\ldots,\int_{[-t,0]}\varphi_N(x+t)d^-\eta(x)\bigg), \\
\sigma(t,\eta) \ &= \ \bar\sigma\bigg(t,\int_{[-t,0]}\varphi_1(x+t)d^-\eta(x),\ldots,\int_{[-t,0]}\varphi_N(x+t)d^-\eta(x)\bigg), \\
F(t,\eta) \ &= \ \bar F\bigg(t,\int_{[-t,0]}\varphi_1(x+t)d^-\eta(x),\ldots,\int_{[-t,0]}\varphi_N(x+t)d^-\eta(x)\bigg), \\
H(\eta) \ &= \ \bar H\bigg(\int_{[-T,0]}\varphi_1(x+T)d^-\eta(x),\ldots,\int_{[-T,0]}\varphi_N(x+T)d^-\eta(x)\bigg),
\end{align*}
where
\begin{itemize}
\item[\textup{(i)}] $\bar b$ and $\bar\sigma$ are continuous functions, with first and second spatial derivatives continuous and satisfying a polynomial growth condition. Moreover, for all $(t,\mathbf x)\in[0,T]\times\R^N$, $\mathbf x_1,\mathbf x_2\in\R^N$,
\begin{align*}
|\bar b(t,\mathbf x)| + |\bar \sigma(t,\mathbf x)| \ &\leq \ \bar M_1(1 + |\mathbf x|), \\
|\bar b(t,\mathbf x_1) - \bar b(t,\mathbf x_2)| + |\bar \sigma(t,\mathbf x_1) - \bar \sigma(t,\mathbf x_2)| \ &\leq \ \bar L_1|\mathbf x_1 - \mathbf x_2|,
\end{align*}
for some positive constants $\bar M_1$ and $\bar L_1$.
\item[\textup{(ii)}] $\bar F$ is continuous and, for all $s\in[0,T]$, the function $\bar F(s,\cdot)$ belongs to $C^2(\R^N)$ and its second order spatial derivatives satisfy a polynomial growth condition uniformly in $s$.
\item[\textup{(iii)}] $\bar H\in C^2(\R^N)$ and its second order spatial derivatives satisfy a polynomial growth condition.
\item[\textup{(iv)}] $\varphi_1,\ldots,\varphi_N\in C^2([0,T])$.
\end{itemize}
Then, the function $\Uc$ given by
\[
\Uc(t,\eta) \ = \ \E\bigg[\int_t^T F(s,\mathbb X_s^{t,\eta})ds + H(\mathbb X_T^{t,\eta})\bigg], \qquad \forall\,(t,\eta)\in[0,T]\times C([-T,0]),
\]
is a strict solution to equation \eqref{KolmEq}.
\end{Theorem}
\textbf{Proof.}
Fix $(t,\eta)\in[0,T]\times C([-T,0])$ and remind that, for any $r\in[t,T]$,
\[
\X_r^{t,\eta}(x) \ = \
\begin{cases}
\eta(r-t+x), \qquad & -T\leq x\leq t-r, \\
X_{r+x}^{t,\eta}, & t-r<x\leq 0.
\end{cases}
\]
Then, can show that, for any $i=1,\ldots,N$ and $r\in[t,T]$, $\P$-a.s. we have
\begin{equation} \label{DetIntStochInt}
\int_{[-r,0]}\varphi_i(x+r)d^-\X_r^{t,\eta}(x) \ = \ \int_{[-t,0]}\varphi_i(x+t)d^-\eta(x) + \int_t^r\varphi_i(u)dX_u^{t,\eta},
\end{equation}
where the left-hand side is intended $\P$-a.s. as a deterministic forward integral. Indeed, the approximation of the left-hand side of \eqref{DetIntStochInt} is given $\P$-a.s. by
\begin{align*}
&\int_{-r-\eps}^0 \varphi_i(x+r) \frac{\X_r^{t,\eta}(x+\eps) - \X_r^{t,\eta}(x)}{\eps} dx\\
&= \ \int_{-r-\eps}^{t-r} \varphi_i(x+r) \frac{\eta(r-t+x+\eps) - \eta(r-t+x)}{\eps} dx + \int_{t-r}^0 \varphi_i(x+r) \frac{X_{r+x+\eps}^{t,\eta} - X_{r+x}^{t,\eta}}{\eps} dx \\
&= \ \int_{-t-\eps}^0 \varphi_i(x+t) \frac{\eta(x+\eps) - \eta(x)}{\eps} dx + \int_t^r \varphi_i(u) \frac{X_{u+\eps}^{t,\eta} - X_u^{t,\eta}}{\eps} du.
\end{align*}
Since $\varphi_i\in C^2([0,T])$, by Proposition \ref{P:BVI}, the left-hand side of the previous chain of equalities goes $\P$-a.s. to
\[
\int_{[-r,0]} \varphi_i(x+r) d^-\X_r^{t,\eta}(x).
\]
Using again Proposition \ref{P:BVI}, we have
\[
\int_{-t-\eps}^0 \varphi_i(x+t) \frac{\eta(x+\eps) - \eta(x)}{\eps} dx \ \overset{\eps\rightarrow0^+}{\longrightarrow} \ \int_{[-r,0]} \varphi_i(x+r) d^-\eta(x).
\]
Finally, taking into account Definition \ref{D:DPFI} and Proposition \ref{Pproperties}(iii), we deduce the following convergence in probability
\[
\int_t^r \varphi_i(u) \frac{X_{u+\eps}^{t,\eta} - X_u^{t,\eta}}{\eps} du \ \overset{\eps\rightarrow0^+}{\longrightarrow} \ \int_t^r \varphi_i(u) dX_u^{t,\eta},
\]
from which \eqref{DetIntStochInt} follows. Therefore,   equation \eqref{SDE} becomes for all $s\in[t,T]$,
\[
\begin{cases}
X_s^{t,\eta} \ = \ \eta(0) + \int_t^s \bar b\bigg(r,\ldots,\int_{[-t,0]}\varphi_i(x+t)d^-\eta(x) + \int_t^r\varphi_i(u)dX_u^{t,\eta},\ldots\bigg) dr \\
\qquad\qquad + \int_t^s \bar\sigma\bigg(r,\ldots,\int_{[-t,0]}\varphi_i(x+t)d^-\eta(x) + \int_t^r\varphi_i(u)dX_u^{t,\eta},\ldots\bigg) dW_r, \hspace{.3cm} t\leq s\leq T, \\
X_s^{t,\eta} \ = \ \eta(s-t), \hspace{8.5cm}-T+t\leq s\leq t.
\end{cases}
\]
For any $(t,\mathbf x)\in[0,T]\times\R^N$, consider the system of equations in $\R^N$ (we denote $\boldsymbol\varphi=(\varphi_1,\ldots,\varphi_N)$)
\begin{equation}
\label{SDE_Y}
\mathbf X_r^{t,\mathbf x} \ = \ \mathbf x + \int_t^r \boldsymbol\varphi(u) \bar b(u,\mathbf X_u^{t,\mathbf x}) du + \int_t^r \boldsymbol\varphi(u) \bar\sigma(u,\mathbf X_u^{t,\mathbf x}) dW_u, \qquad \forall\,r\in[t,T].
\end{equation}
Under assumption (i) on $\bar b$ and $\bar\sigma$, it is well-known (see, e.g., Theorem 1.1, Chapter 5, in \cite{friedman75vol1}) that there exists a  unique (up to indistinguishability) continuous process $(X^{t,\mathbf x,1},\ldots,X^{t,\mathbf x,N})$ $=\mathbf X^{t,\mathbf x}=(\mathbf X_r^{t,\mathbf x})_{r\in[t,T]}$ solution to \eqref{SDE_Y}. Notice that, when $\mathbf x=(x_1,\ldots,x_N)\in\R^N$ is given by
\[
x_i \ = \ \int_{[-t,0]}\varphi_i(x+t)d^-\eta(x),
\]
then, by uniqueness of \eqref{SDE_Y} and \eqref{DetIntStochInt},
\begin{equation}
\label{Proof}
X_r^{t,\mathbf x,i} \ = \ \int_{[-r,0]}\varphi_i(x+r)d^-\X_r^{t,\eta}(x) \ = \ \int_{[-t,0]}\varphi_i(x+t)d^-\eta(x) + \int_t^r\varphi_i(u)dX_u^{t,\eta}.
\end{equation}
As a consequence, we obtain
\begin{align*}
\Uc(t,\eta) \ &= \ \E\bigg[\int_t^T F(s,\mathbb X_s^{t,\eta})ds + H(\mathbb X_T^{t,\eta})\bigg] \\
&= \ \E\bigg[\int_t^T \bar F\bigg(s,\ldots,\int_{[-t,0]}\varphi_i(x+t)d^-\eta(x) + \int_t^s\varphi_i(u)dX_u^{t,\eta},\ldots\bigg) ds \\
&\quad \ + \bar H\bigg(\ldots,\int_{[-t,0]}\varphi_i(x+t)d^-\eta(x) + \int_t^T\varphi_i(u)dX_u^{t,\eta},\ldots\bigg)\bigg] \\
&= \ \Psi\bigg(t,\int_{[-t,0]}\varphi_1(x+t)d^-\eta(x),\ldots,\int_{[-t,0]}\varphi_N(x+t)d^-\eta(x)\bigg),
\end{align*}
where
\[
\Psi(t,\mathbf x) \ = \ \E\bigg[\int_t^T \bar F\big(s,\mathbf X_s^{t,\mathbf x}\big) ds + \bar H\big(\mathbf X_T^{t,\mathbf x}\big)\bigg], \qquad \forall\,(t,\mathbf x)\in[0,T]\times\R^N.
\]
Notice that $\Psi\in C([0,T]\times\R^N)$, as a consequence of the continuous dependence of $\mathbf X_s^{t,\mathbf x}$ on $(t,\mathbf x)$ and of the standard estimate $\sup_{t\in[0,T],\,|\mathbf x|\leq R}\E[\sup_{s\in[t,T]}|\mathbf X_s^{t,\mathbf x}|^p]<\infty$, for any $p\geq1$ and $R>0$. If $\bar F\equiv0$, it follows from Theorem 6.1, Chapter 5, in \cite{friedman75vol1}, that $\Psi\in C^{1,2}([0,T]\times\R^N)$ and satisfies the following backward parabolic equation:
\[
\begin{cases}
\partial_t \Psi(t,\mathbf x) + \sum_{i=1}^N \varphi_i(t)\bar b(t,\mathbf x)D_{x_i}\Psi(t,\mathbf x) \\
+ \frac{1}{2}\sum_{i,j=1}^N \varphi_i(t)\varphi_j(t) \bar\sigma^2(t,\mathbf x) D_{x_ix_j}^2 \Psi(t,\mathbf x) = 0, \qquad\qquad &\forall\,(t,\mathbf x)\in[0,T[\times\R^N, \\
\Psi(T,\mathbf x) = \bar H(\mathbf x), &\forall\,\mathbf x\in\R^N.
\end{cases}
\]
When $\bar F\not\equiv0$, the result is still true and the proof is based on Duhamel's principle. More precisely, we fix $s\in[0,T]$ and consider the following equation:
\begin{equation}
\label{PDE_Psi^s}
\begin{cases}
\partial_t \Psi^s(t,\mathbf x) + \sum_{i=1}^N \varphi_i(t)\bar b(t,\mathbf x)D_{x_i}\Psi^s(t,\mathbf x) \\
+ \frac{1}{2}\sum_{i,j=1}^N \varphi_i(t)\varphi_j(t) \bar\sigma^2(t,\mathbf x) D_{x_ix_j}^2 \Psi^s(t,\mathbf x) = 0, \qquad\qquad &\forall\,(t,\mathbf x)\in[0,s[\times\R^N, \\
\Psi^s(s,\mathbf x) = \bar F(s,\mathbf x), &\forall\,\mathbf x\in\R^N.
\end{cases}
\end{equation}
Let $\Delta:=\{(s,t)\in[0,T]\times[0,T]\colon0\leq t\leq s\}$ and define the map
\begin{align*}
\Psi^\cdot(\cdot,\cdot)\colon\Delta\times\R^N&\rightarrow\R \\
(s,t,\mathbf x)&\mapsto\Psi^s(t,\mathbf x) \ = \ \E\big[\bar F\big(s,\mathbf X_s^{t,\mathbf x}\big)\big], \qquad \forall\,(t,\mathbf x)\in[0,s]\times\R^N.
\end{align*}
Notice that $\Psi^\cdot(\cdot,\cdot)\in C(\Delta\times\R^N)$, as a consequence of the continuous dependence of $\mathbf X_s^{t,\mathbf x}$ on $(s,t,\mathbf x)$ and of the standard estimate $\sup_{t\in[0,T],\,|\mathbf x|\leq R}\E[\sup_{r\in[t,T]}|\mathbf X_r^{t,\mathbf x}|^p]<\infty$, for any $p\geq1$ and $R>0$. Using again Theorem 6.1, Chapter 5, in \cite{friedman75vol1}, we see that, for any $s\in[0,T]$, the function $\Psi^s$ belongs to $C^{1,2}([0,s]\times\R^N)$ and satisfies equation \eqref{PDE_Psi^s}. Moreover, from the proof of Theorem 5.5, Chapter 5, in \cite{friedman75vol1}, we see that, for any $i=1,\ldots,N$, the map $D_{x_i}\Psi^\cdot(\cdot,\cdot)\colon\Delta\times\R^N\rightarrow\R$ is continuous. As a consequence, it follows from \eqref{PDE_Psi^s} that the map $\partial_t\Psi^\cdot(\cdot,\cdot)\colon\Delta\times\R^N\rightarrow\R$ is also continuous. Then, by direct calculation, we see that $\Psi$, which can be written as
\[
\Psi(t,\mathbf x) \ = \ \E\bigg[\int_t^T \bar F\big(s,\mathbf X_s^{t,\mathbf x}\big) ds + \bar H\big(\mathbf X_T^{t,\mathbf x}\big)\bigg] \ = \ \int_t^T \Psi^s(t,\mathbf x) ds + \E\big[\bar H\big(\mathbf X_T^{t,\mathbf x}\big)\big],
\]
is a classical solution to the backward parabolic PDE
\begin{equation}
\label{PDE_Psi}
\begin{cases}
\partial_t \Psi(t,\mathbf x) + \sum_{i=1}^N \varphi_i(t)\bar b(t,\mathbf x)D_{x_i}\Psi(t,\mathbf x) \\
+ \frac{1}{2}\sum_{i,j=1}^N \varphi_i(t)\varphi_j(t) \bar\sigma^2(t,\mathbf x) D_{x_ix_j}^2 \Psi(t,\mathbf x) + \bar F(t,\mathbf x) = 0, &\quad\forall\,(t,\mathbf x)\in[0,T[\times\R^N, \\
\Psi(T,\mathbf x) = \bar H(\mathbf x), &\quad\forall\,\mathbf x\in\R^N.
\end{cases}
\end{equation}
Finally, we derive formulae for the derivatives of $\Uc$, expressed in terms of the derivatives of $\Psi$. We begin noting that, taking into account Proposition \ref{P:BVI}, we have
\[
\int_{[-t,0]}\varphi_i(x+t)d^-\eta(x) \ = \ \eta(0)\varphi_i(t) - \int_{-t}^0 \eta(x)\dot{\varphi}_i(x+t)dx, \qquad \forall\,\eta\in C([-T,0]).
\]
This in turn implies that $\Uc$ admits a continuous extension (necessarily 
unique) $u\colon\mathscr C([-T,0])$ $\rightarrow\R$ given by
\[
u(t,\eta) \ = \ \Psi\bigg(t,\int_{[-t,0]}\varphi_1(x+t)d^-\eta(x),\ldots,\int_{[-t,0]}\varphi_N(x+t)d^-\eta(x)\bigg),
\]
for all $(t,\eta)\in[0,T]\times\mathscr C([-T,0])$. We also define the map $\tilde u\colon[0,T]\times\mathscr C([-T,0[)\times\R\rightarrow\R$ as in \eqref{E:tildeu}:
\[
\tilde u(t,\gamma,a) \ = \ u(t,\gamma1_{[-T,0[}+a1_{\{0\}}) \ = \ \Psi\bigg(t,\ldots,a\varphi_i(t) - \int_{-t}^0 \gamma(x)\dot{\varphi}_i(x+t)dx,\ldots\bigg),
\]
for all $(t,\gamma,a)\in[0,T]\times\mathscr C([-T,0[)\times\R$. Let us evaluate the time derivative $\partial_t\Uc(t,\eta)$, for a given $(t,\eta)\in[0,T[\times C([-T,0])$:
\begin{align*}
\partial_t\Uc(t,\eta) \ &= \ \partial_t\Psi\bigg(t,\int_{[-t,0]}\varphi_1(x+t)d^-\eta(x),\ldots,\int_{[-t,0]}\varphi_N(x+t)d^-\eta(x)\bigg) \\
&\quad \ + \sum_{i=1}^N D_{x_i}\Psi\bigg(t,\ldots,\int_{[-t,0]}\varphi_i(x+t)d^-\eta(x),\ldots\bigg)\partial_t\bigg(\int_{[-t,0]}\varphi_i(x+t)d^-\eta(x)\bigg).
\end{align*}
Notice that
\begin{align*}
\partial_t\bigg(\int_{[-t,0]}\varphi_i(x+t)d^-\eta(x)\bigg) \ &= \ \partial_t\bigg(\eta(0)\varphi(t) - \int_{-t}^0\eta(x)\dot\varphi_i(x+t)dx\bigg) \\
&= \ \eta(0)\dot\varphi(t) - \eta(-t)\dot\varphi_i(0^+) - \int_{-t}^0 \eta(x)\ddot\varphi_i(x+t) dx.
\end{align*}
Let us proceed with the horizontal derivative.  We have
\begin{align*}
&D^H\Uc(t,\eta) \ = \ D^H u(t,\eta) \ = \ D^H\tilde u(t,\eta_{|[-T,0[},\eta(0)) \\
&= \ \lim_{\eps\rightarrow0^+} \frac{\tilde u(t,\eta_{|[-T,0[}(\cdot),\eta(0)) - \tilde u(t,\eta_{|[-T,0[}(\cdot-\eps),\eta(0))}{\eps} \\
&= \ \lim_{\eps\rightarrow0^+} \bigg(\frac{1}{\eps}\Psi\left(t,\ldots,\eta(0)\varphi_i(t) - \int_{-t}^0 \eta(x)\dot{\varphi}_i(x+t)dx,\ldots\right) \\
&\quad \ - \frac{1}{\eps}\Psi\left(t,\ldots,\eta(0)\varphi_i(t) - \int_{-t}^0 \eta(x-\eps)\dot{\varphi}_i(x+t)dx,\ldots\right)\bigg).
\end{align*}
From the fundamental theorem of calculus, we obtain
\begin{align*}
&\frac{1}{\eps}\Psi\left(t,\ldots,\eta(0)\varphi_i(t) - \int_{-t}^0 \eta(x)\dot{\varphi}_i(x+t)dx,\ldots\right) \\
&- \frac{1}{\eps}\Psi\left(t,\ldots,\eta(0)\varphi_i(t) - \int_{-t}^0 \eta(x-\eps)\dot{\varphi}_i(x+t)dx,\ldots\right) \\
&= \ \frac{1}{\eps}\int_0^\eps \sum_{i=1}^N D_{x_i}\Psi\left(t,\ldots,\eta(0)\varphi_i(t) - \int_{-t}^0 \eta(x-y)\dot{\varphi}_i(x+t)dx,\ldots\right)\partial_y\bigg(\eta(0)\varphi_i(t) \\
&\quad \ - \int_{-t}^0 \eta(x-y)\dot{\varphi}_i(x+t)dx\bigg) dy.
\end{align*}
Notice that
\begin{align*}
&\partial_y\bigg(\eta(0)\varphi_i(t) - \int_{-t}^0 \eta(x-y)\dot{\varphi}_i(x+t)dx\bigg) \ = \ -\partial_y\bigg(\int_{-t-y}^{-y}\eta(x)\dot{\varphi}_i(x+y+t)dx\bigg) \\
&= \ - \bigg(\eta(-y)\dot{\varphi}_i(t) - \eta(-t-y)\dot{\varphi}_i(0^+) + \int_{-t-y}^{-y}\eta(x)\ddot{\varphi}_i(x+y+t)dx\bigg).
\end{align*}
Therefore
\begin{align*}
D^H\Uc(t,\eta) \ &= \ -\lim_{\eps\rightarrow0^+} \frac{1}{\eps}\int_0^\eps \sum_{i=1}^N D_{x_i}\Psi\bigg(t,\ldots,\eta(0)\varphi_i(t) - \int_{-t}^0\!\! \eta(x-y)\dot{\varphi}_i(x+t)dx,\ldots\bigg)\bigg(\eta(-y)\dot{\varphi}_i(t) \\
&\quad \ - \eta(-t-y)\dot{\varphi}_i(0^+) + \int_{-t-y}^{-y}\eta(x)\ddot{\varphi}_i(x+y+t)dx\bigg) dy \\
&=\ - \sum_{i=1}^N D_{x_i}\Psi\bigg(t,\ldots,\eta(0)\varphi_i(t) - \int_{-t}^0 \eta(x)\dot{\varphi}_i(x+t)dx,\ldots\bigg)\bigg(\eta(0)\dot\varphi(t) - \eta(-t)\dot\varphi_i(0^+) \\
&\quad \ - \int_{-t}^0 \eta(x)\ddot\varphi_i(x+t) dx\bigg).
\end{align*}
Finally, concerning the vertical derivative we have
\begin{align*}
D^V\Uc(t,\eta) \ = \ D^Vu(t,\eta) \ &= \ \partial_a\tilde u(t,\eta1_{[-T,0[}+\eta(0)1_{\{0\}}) \\
&= \ \sum_{i=1}^N D_{x_i}\Psi\bigg(t,\int_{[-t,0]}\varphi_1(x+t)d^-\eta(x),\ldots\bigg)\varphi_i(t)
\end{align*}
and
\begin{align*}
D^{VV}\Uc(t,\eta) \ = \ D^{VV}u(t,\eta) \ &= \ \partial_{aa}^2\tilde u(t,\eta1_{[-T,0[}+\eta(0)1_{\{0\}}) \\
&= \ \sum_{i,j=1}^N D_{x_ix_j}^2\Psi\bigg(t,\int_{[-t,0]}\varphi_1(x+t)d^-\eta(x),\ldots\bigg)\varphi_i(t)\varphi_j(t).
\end{align*}
From the regularity of $\Psi$, we see that $\Uc\in C^{1,2}(([0,T[\times\textup{past})\times\textup{present}))\cap C([0,T]\times C([-T,0]))$. Furthermore, as $\Psi$ is a solution to \eqref{PDE_Psi}, it follows that $\Uc$ solves equation \eqref{KolmEq}.
\ep

\begin{Theorem}
\label{T:ExistenceStrictPardouxPeng}
Suppose that there exists $N\in\N\backslash\{0\}$ such that, for all $(t,\eta,y,z)\in[0,T]\times C([-T,0])$ $\times\R\times\R$,
\begin{align*}
b(t,\eta) \ &= \ \bar b\bigg(\int_{[-t,0]}\varphi_1(x+t)d^-\eta(x),\ldots,\int_{[-t,0]}\varphi_N(x+t)d^-\eta(x)\bigg), \\
\sigma(t,\eta) \ &= \ \bar\sigma\bigg(\int_{[-t,0]}\varphi_1(x+t)d^-\eta(x),\ldots,\int_{[-t,0]}\varphi_N(x+t)d^-\eta(x)\bigg), \\
F(t,\eta,y,z) \ &= \ \bar F\bigg(t,\int_{[-t,0]}\varphi_1(x+t)d^-\eta(x),\ldots,\int_{[-t,0]}\varphi_N(x+t)d^-\eta(x),y,z\bigg), \\
H(\eta) \ &= \ \bar H\bigg(\int_{[-T,0]}\varphi_1(x+T)d^-\eta(x),\ldots,\int_{[-T,0]}\varphi_N(x+T)d^-\eta(x)\bigg),
\end{align*}
where
\begin{itemize}
\item[\textup{(i)}] $\bar b$, $\bar\sigma$, $\bar F$, $\bar H$ are continuous functions satisfying, for some positive constants $C$ and $m$,
\begin{align*}
|\bar b(\mathbf x)-\bar b(\mathbf x')| + |\bar\sigma(\mathbf x)-\bar\sigma(\mathbf x')| \ &\leq \ C|\mathbf x-\mathbf x'|, \\
|\bar F(t,\mathbf x,y,z)-\bar F(t,\mathbf x,y',z')| \ &\leq \ C\big(|y-y'| + |z-z'|\big), \\
|\bar F(t,\mathbf x,0,0)| + |\bar H(\mathbf x)| \ &\leq \ C\big(1 + |\mathbf x|^m\big),
\end{align*}
for all $t\in[0,T]$, $\mathbf x,\mathbf x'\in\R^N$, $y,y'\in\R$, and $z,z'\in\R$.
\item[\textup{(ii)}] $\bar b$ and $\bar\sigma$ are of class $C^3$ with partial derivatives from order $1$ up to order $3$ bounded.
\item[\textup{(iii)}] For all $t\in[0,T]$, $\bar F(t,\cdot,\cdot,\cdot)\in C^3(\R^N)$ and moreover the following.
\begin{enumerate}
\item[\textup{(a)}] $\bar F(t,\cdot,0,0)\in C^3(\R^N)$ and its third order partial derivatives satisfy a polynomial growth condition uniformly in $t$.
\item[\textup{(b)}] $D_y\bar F$, $D_z\bar F$ are bounded on $[0,T]\times\R^N\times\R\times\R$, as well as their derivatives of order one and second with respect to $x_1,\ldots,x_N,y,z$.
\end{enumerate}
\item[\textup{(iv)}] $\bar H\in C^3(\R^N)$ and its third order partial derivatives satisfy a polynomial growth condition.
\item[\textup{(v)}] $\varphi_1,\ldots,\varphi_N\in C^2([0,T])$.
\end{itemize}
Then, the map $\Uc$ given by
\[
\Uc(t,\eta) \ = \ Y_t^{t,\eta}, \qquad \forall\,(t,\eta)\in[0,T]\times C([-T,0]),
\]
where $(Y_s^{t,\eta},Z_s^{t,\eta})_{s\in[t,T]}\in\S^2(t,T)\times\H^2(t,T)$ is the unique solution to \eqref{BSDE}, is a strict solution to equation \eqref{KolmEq}.
\end{Theorem}
\textbf{Proof.}
The proof can be done proceeding as in the proof of Theorem \ref{T:ExistenceStrict}, and it is based on Theorem 3.2 in \cite{pardoux_peng92} instead of Theorem 6.1, Chapter 5, in \cite{friedman75vol1}. More precisely, adopting the same notations as in the proof of Theorem \ref{T:ExistenceStrict}, set $\boldsymbol\varphi=(\varphi_1,\ldots,\varphi_N)$ and consider, for any $(t,\mathbf x)\in[0,T]\times\R^N$, the forward-backward system of stochastic differential equations:
\begin{equation}
\label{FBSDE}
\begin{cases}
\mathbf X_s^{t,\mathbf x} \ = \ \mathbf x + \int_t^s \boldsymbol\varphi(r) \bar b(r,\mathbf X_r^{t,\mathbf x}) dr + \int_t^s \boldsymbol\varphi(r) \bar\sigma(r,\mathbf X_r^{t,\mathbf x}) dW_r, &\quad s\in[t,T], \\
Y_s^{t,\mathbf x} \ = \ \bar H(\mathbf X_T^{t,\mathbf x}) + \int_s^T \bar F(r,\mathbf X_r^{t,\mathbf x},Y_r^{t,\mathbf x},Z_r^{t,\mathbf x}) dr - \int_s^T Z_r^{t,\mathbf x} dW_r, &\quad s\in[t,T].
\end{cases}
\end{equation}
Under assumption (i) on $\bar b$ and $\bar\sigma$, we see that that there exists a  unique (up to indistinguishability) continuous process $(X^{t,\mathbf x,1},\ldots,X^{t,\mathbf x,N})=\mathbf X^{t,\mathbf x}=(\mathbf X_s^{t,\mathbf x})_{s\in[t,T]}$ solution to the forward equation in \eqref{FBSDE}, see, e.g., Theorem 1.1, Chapter 5, in \cite{friedman75vol1}. Moreover, from Theorem 4.1 in \cite{parpen90} it follows that, under assumption (i) on $\bar F$ and $\bar H$, there exists a unique solution $(Y^{t,\mathbf x},Z^{t,\mathbf x})\in\S^2(t,T)\times\H^2(t,T)$ to the backward equation in \eqref{FBSDE}. Now, fix $\eta\in C([-T,0])$ and define $\mathbf x=(x_1,\ldots,x_N)\in\R^N$ as
\begin{equation}
\label{x=eta}
x_i \ = \ \int_{[-t,0]}\varphi_i(x+t)d^-\eta(x).
\end{equation}
Then, similarly as for \eqref{Proof},
\[
X_r^{t,\mathbf x,i} \ = \ \int_{[-t,0]}\varphi_i(x+t)d^-\eta(x) + \int_t^r\varphi_i(u)dX_u^{t,\eta} \ = \ \int_{[-r,0]}\varphi_i(x+r)d^-\X_r^{t,\eta}(x)
\]
and therefore
\[
Y_s^{t,\mathbf x} \ = \ H(\mathbb X_T^{t,\eta}) + \int_s^T F(r,\mathbb X_r^{t,\eta},Y_r^{t,\mathbf x},Z_r^{t,\mathbf x}) dr - \int_s^T Z_r^{t,\mathbf x} dW_r, \qquad t\leq s\leq T.
\]
Since $(Y^{t,\eta},Z^{t,\eta})$ also solves the above backward equation, from the uniqueness to the BSDE it follows that $Y^{t,\mathbf x}=Y^{t,\eta}$ in $\S^2(t,T)$ and $Z^{t,\mathbf x}=Z^{t,\eta}$ and $\H^2(t,T)$, whenever $\mathbf x=(x_i)_{i=1,\ldots,N}$ is given by \eqref{x=eta}. Now, from the definition of $\Uc$ and the equation \eqref{BSDE} satisfied by $(Y^{t,\eta},Z^{t,\eta})$, we have (with $\eta$ and $\mathbf x$ related by \eqref{x=eta})
\begin{align*}
\Uc(t,\eta) \ = \ Y_t^{t,\eta} \ &= \ \E\bigg[\int_t^T F(s,\mathbb X_s^{t,\eta},Y_s^{t,\eta},Z_s^{t,\eta})ds + H(\mathbb X_T^{t,\eta})\bigg] \\
&= \ \E\bigg[\int_t^T \bar F(s,\mathbf X_s^{t,\mathbf x},Y_s^{t,\mathbf x},Z_s^{t,\mathbf x}) ds + \bar H(\mathbf X_T^{t,\mathbf x})\bigg].
\end{align*}
Then, we define
\[
\Psi(t,\mathbf x) \ = \ \E\bigg[\int_t^T \bar F(s,\mathbf X_s^{t,\mathbf x},Y_s^{t,\mathbf x},Z_s^{t,\mathbf x}) ds + \bar H(\mathbf X_T^{t,\mathbf x})\bigg], \qquad \forall\,(t,\mathbf x)\in[0,T]\times\R^N.
\]
It follows from Theorem 3.2 in \cite{pardoux_peng92} that $\Psi\in C^{1,2}([0,T]\times\R^N)$ and satisfies the following backward semilinear parabolic equation
\[
\begin{cases}
\partial_t \Psi(t,\mathbf x) + \sum_{i=1}^N \varphi_i(t)\bar b(t,\mathbf x)D_{x_i}\Psi(t,\mathbf x) + \frac{1}{2}\sum_{i,j=1}^N \varphi_i(t)\varphi_j(t) \bar\sigma^2(t,\mathbf x) D_{x_ix_j}^2 \Psi(t,\mathbf x) \\
+ \bar F\big(t,\mathbf x,\Psi(t,\mathbf x),\sum_{i=1}^N\bar\sigma(t,\mathbf x)\varphi_i(t)D_{x_i}\Psi(t,\mathbf x)\big) = 0, \hspace{3.7cm}\forall\,(t,\mathbf x)\in[0,T[\times\R^N, \\
\Psi(T,\mathbf x) = \bar H(\mathbf x), \hspace{9.35cm}\forall\,\mathbf x\in\R^N.
\end{cases}
\]
Finally, the claim follows expressing the derivatives of $\Uc$ in terms of the derivatives of $\Psi$ as in the proof of Theorem \ref{T:ExistenceStrict}.
\ep

\vspace{5mm}

\noindent\textbf{Acknowledgements.} Part of the paper was done during the visit
of the second named author at the Centre for Advanced Study (CAS) at the Norwegian Academy of Science and Letters in Oslo.
That  author also benefited partially from the
support of the ``FMJH Program Gaspard Monge in optimization and operation
research'' (Project 2014-1607H).

\small
\bibliographystyle{plain}
\bibliography{biblio}

\def\polhk#1{\setbox0=\hbox{#1}{\ooalign{\hidewidth
  \lower1.5ex\hbox{`}\hidewidth\crcr\unhbox0}}}
  \def\polhk#1{\setbox0=\hbox{#1}{\ooalign{\hidewidth
  \lower1.5ex\hbox{`}\hidewidth\crcr\unhbox0}}} \def\cprime{$'$}
  \def\polhk#1{\setbox0=\hbox{#1}{\ooalign{\hidewidth
  \lower1.5ex\hbox{`}\hidewidth\crcr\unhbox0}}}
\begin{thebibliography}{10}

\bibitem{buckdahn_ma_zhang13}
R.~Buckdahn, J.~Ma, and J.~Zhang.
\newblock {\em \emph{Pathwise Taylor Expansions for Random Fields on Multiple
  Dimensional Paths}}.
\newblock \emph{Preprint} arXiv:1310.0517, 2013.

\bibitem{cho}
A.~Chojnowska-Michalik.
\newblock Representation theorem for general stochastic delay equations.
\newblock {\em Bull. Acad. Polon. Sci. S\'er. Sci. Math. Astronom. Phys.},
  26(7):635--642, 1978.

\bibitem{contfournie10}
R.~Cont and D.-A. Fourni{\'e}.
\newblock Change of variable formulas for non-anticipative functionals on path
  space.
\newblock {\em J. Funct. Anal.}, 259(4):1043--1072, 2010.

\bibitem{contfournie}
R.~Cont and D.-A. Fourni{\'e}.
\newblock A functional extension of the {I}t\^o formula.
\newblock {\em C. R. Math. Acad. Sci. Paris}, 348(1-2):57--61, 2010.

\bibitem{contfournie13}
R.~Cont and D.-A. Fourni{\'e}.
\newblock Functional {I}t\^o calculus and stochastic integral representation of
  martingales.
\newblock {\em Ann. Probab.}, 41(1):109--133, 2013.

\bibitem{cosso_russo15b}
A.~Cosso and F.~Russo.
\newblock {\em \emph{Strong-viscosity solutions: semilinear parabolic {PDE}s
  and path-dependent {PDE}s}}.
\newblock \emph{Preprint}, 2015.

\bibitem{digirfabbrirusso13}
C.~Di~Girolami, G.~Fabbri, and F.~Russo.
\newblock The covariation for {B}anach space valued processes and applications.
\newblock {\em Metrika}, 77(1):51--104, 2014.

\bibitem{DGR}
C.~Di~Girolami and F.~Russo.
\newblock Infinite dimensional stochastic calculus via regularization and
  applications.
\newblock {\em Preprint \textup{HAL-INRIA, inria-00473947 version 1}},
  (Unpublished), 2010.

\bibitem{DGRnote}
C.~Di~Girolami and F.~Russo.
\newblock Clark-{O}cone type formula for non-semimartingales with finite
  quadratic variation.
\newblock {\em C. R. Math. Acad. Sci. Paris}, 349(3-4):209--214, 2011.

\bibitem{DGR2}
C.~Di~Girolami and F.~Russo.
\newblock Generalized covariation and extended {F}ukushima decomposition for
  {B}anach space-valued processes. {A}pplications to windows of {D}irichlet
  processes.
\newblock {\em Infin. Dimens. Anal. Quantum Probab. Relat. Top.},
  15(2):1250007, 50, 2012.

\bibitem{digirrusso12}
C.~Di~Girolami and F.~Russo.
\newblock Generalized covariation for {B}anach space valued processes, {I}t\^o
  formula and applications.
\newblock {\em Osaka J. Math.}, 51(3):729--783, 2014.

\bibitem{dupire}
B.~Dupire.
\newblock {\em \emph{Functional {I}t\^o calculus}}.
\newblock \emph{Portfolio Research Paper}, Bloomberg, 2009.

\bibitem{flandoli_zanco13}
F.~Flandoli and G.~Zanco.
\newblock {\em \emph{An infinite-dimensional approach to path-dependent
  Kolmogorov's equations}}.
\newblock \emph{Preprint} arXiv:1312.6165, 2013.

\bibitem{follmer}
H.~F{\"o}llmer.
\newblock Calcul d'{I}t\^o sans probabilit\'es.
\newblock In {\em Seminar on {P}robability, {XV} ({U}niv. {S}trasbourg,
  {S}trasbourg, 1979/1980) ({F}rench)}, volume 850 of {\em Lecture Notes in
  Math.}, pages 143--150. Springer, Berlin, 1981.

\bibitem{friedman75vol1}
A.~Friedman.
\newblock {\em Stochastic differential equations and applications. {V}ol. 1}.
\newblock Academic Press, New York, 1975.
\newblock Probability and Mathematical Statistics, Vol. 28.

\bibitem{gradinaru_nourdin03}
M.~Gradinaru and I.~Nourdin.
\newblock Approximation at first and second order of {$m$}-order integrals of
  the fractional {B}rownian motion and of certain semimartingales.
\newblock {\em Electron. J. Probab.}, 8:no. 18, 26 pp. (electronic), 2003.

\bibitem{jacod79}
J.~Jacod.
\newblock {\em Calcul stochastique et probl\`emes de martingales}, volume 714
  of {\em Lecture Notes in Mathematics}.
\newblock Springer, Berlin, 1979.

\bibitem{leao_ohashi_simas14}
D.~Le\~ao, A.~Ohashi, and A.~B. Simas.
\newblock {\em \emph{Weak functional {I}t\^o calculus and applications}}.
\newblock \emph{Preprint} arXiv:1408.1423v2, 2014.

\bibitem{parpen90}
{\'E}.~Pardoux and S.~Peng.
\newblock Adapted solution of a backward stochastic differential equation.
\newblock {\em Systems Control Lett.}, 14(1):55--61, 1990.

\bibitem{pardoux_peng92}
{\'E}.~Pardoux and S.~Peng.
\newblock Backward stochastic differential equations and quasilinear parabolic
  partial differential equations.
\newblock In {\em Stochastic partial differential equations and their
  applications ({C}harlotte, {NC}, 1991)}, volume 176 of {\em Lecture Notes in
  Control and Inform. Sci.}, pages 200--217. Springer, Berlin, 1992.

\bibitem{pazy83}
A.~Pazy.
\newblock {\em Semigroups of linear operators and applications to partial
  differential equations}, volume~44 of {\em Applied Mathematical Sciences}.
\newblock Springer-Verlag, New York, 1983.

\bibitem{rudin}
W.~Rudin.
\newblock {\em Real and complex analysis}.
\newblock McGraw-Hill Book Co., New York, third edition, 1987.

\bibitem{russovallois91}
F.~Russo and P.~Vallois.
\newblock Int\'egrales progressive, r\'etrograde et sym\'etrique de processus
  non adapt\'es.
\newblock {\em C. R. Acad. Sci. Paris S\'er. I Math.}, 312(8):615--618, 1991.

\bibitem{russovallois93}
F.~Russo and P.~Vallois.
\newblock Forward, backward and symmetric stochastic integration.
\newblock {\em Probab. Theory Related Fields}, 97(3):403--421, 1993.

\bibitem{russovallois95}
F.~Russo and P.~Vallois.
\newblock The generalized covariation process and {I}t\^o formula.
\newblock {\em Stochastic Process. Appl.}, 59(1):81--104, 1995.

\bibitem{russovallois07}
F.~Russo and P.~Vallois.
\newblock Elements of stochastic calculus via regularization.
\newblock In {\em S\'eminaire de {P}robabilit\'es {XL}}, volume 1899 of {\em
  Lecture Notes in Math.}, pages 147--185. Springer, Berlin, 2007.

\bibitem{ryan02}
R.~A. Ryan.
\newblock {\em Introduction to tensor products of {B}anach spaces}.
\newblock Springer Monographs in Mathematics. Springer-Verlag London Ltd.,
  London, 2002.

\end{thebibliography}

\end{document}